\documentclass[11pt]{article}
\usepackage{amssymb}
\usepackage{colortbl}
\usepackage{color}

\newtheorem{theorem}{Theorem}[section]

\newtheorem{corollary}[theorem]{Corollary}
\newtheorem{definition}[theorem]{Definition}

\newtheorem{lemma}[theorem]{Lemma}

\begin{document}

\renewcommand{\theequation}{\arabic{section}.\arabic{equation}}

\bigskip

\begin{center}
{\Large An operatorial approach of the well-posedness of an algebraic
Riccati equation}

\bigskip

Gabriela Marinoschi

\medskip

\textquotedblleft Gheorghe Mihoc-Caius Iacob\textquotedblright\ Institute of
Mathematical Statistics and

Applied Mathematics of the Romanian Academy,

Calea 13 Septembrie 13, Bucharest, Romania

E-mail: gabriela.marinoschi@acad.ro

\bigskip
\end{center}

\noindent \noindent Abstract. Finding the state feedback control in an $%
H^{\infty }$-optimal control problem involves a challenging approach of the
associated algebraic Riccati equation of the generic form $A^{\ast
}P+PA+P\Gamma P=F$. In view of this objective, we explore in this paper the
existence of the solution to this algebraic Riccati equation by a direct
operatorial approach in the space of Hilbert-Schmidt operators. The proofs
are provided, under certain assumptions on the operators $\Gamma $ and $F,$
for the cases with $A$ coercive and $A\geq 0,$ respectively. They develop a
constructive approach, possibly indicating a method for finding the
numerical solution. Next, relying on the existence of the solution to the
Riccati equation, we provide then a result concerning the associated $%
H^{\infty }$-optimal control problem. An example regarding the application
of the existence proof for the solution to the Riccati equation is given for
a parabolic equation with a singular potential of Hardy type.

\medskip


\medskip


\section{Introduction}

\setcounter{equation}{0}

Riccati equation plays an essential role in the robust stabilization of
partial differential equations which are corrupted by perturbations. Namely,
the feedback control found to robustly stabilize a perturbed system of PDEs,
via the $H^{\infty }$-optimal control technique, is constructed by means of
the solution to an associated algebraic Riccati equation.

We shall briefly describe the problem only for justifying the interest in
studying the existence of the solution to a Riccati equation in the infinite
dimensional case. Let $H,$ $W,$ $Z,$ $U$ be Hilbert spaces identified with
their duals, let $B_{1}:W\rightarrow H,$ $B_{2}:U\rightarrow H,$ $%
C_{1}:H\rightarrow Z,$ $D_{1}:U\rightarrow Z$ be linear and bounded
operators, let $A$ be a linear operator defined on $D(A)\subset H\rightarrow
H$ such that $-A$ is a generator of a $C_{0}$-semigroup on $H.$ We consider
the system 
\begin{equation}
y^{\prime }(t)+Ay(t)=B_{1}w(t)+B_{2}u(t),\mbox{ }t>0,  \label{1}
\end{equation}%
\begin{equation}
z(t)=\Phi (y(t),u(t))\mbox{ }t>0,  \label{2}
\end{equation}%
\[
y(0)=0, 
\]%
where $y$ is the state, $w$ is a unknown perturbation, $u$ is the control to
be found and $z$ is an output. Here, we consider a general case for the
output form, where 
\[
\Phi :L^{2}(\mathbb{R}_{+};H)\times L^{2}(\mathbb{R}_{+};U)\rightarrow L^{2}(%
\mathbb{R}_{+};Z),\mbox{ }\mathbb{R}_{+}=(0,\infty ), 
\]%
such that 
\begin{equation}
\left\Vert \Phi (y(t),u(t))\right\Vert _{Z}^{2}\leq \left\Vert
C_{1}y(t)\right\Vert _{Z}^{2}+\left\Vert u(t)\right\Vert _{U}^{2}.
\label{2-1}
\end{equation}%
{%
\color{black}%
}This form covers the particular choices $z(t)=C_{1}y(t)+D_{1}u(t)$ (with
the conditions $D_{1}^{\ast }C_{1}=0$ and $D_{1}D_{1}^{\ast }=1,$ used in 
\cite{GM-23-Hinf}), as well as the generally used output form $%
z(t)=C_{1}y(t).$

The robust stabilization of this system, via the $H^{\infty }$-optimal
control technique, means to design a feedback control $u=Fy$ that stabilizes
exponentially the system and helps it to achieve a robust performance,
ensuring the fact that the perturbation may not have an extreme influence on
the output. More precisely, the problem is to determine $F:H\rightarrow U$
such that the operator $-A+B_{2}F$ generates an exponentially stable
semigroup on $H,$ namely that 
\begin{equation}
\left\Vert e^{(-A+B_{2}F)t}y\right\Vert _{H}\leq Ce^{-\alpha t}\left\Vert
y\right\Vert _{H}\mbox{ for all }y\in H,\mbox{ with }C>0,\mbox{ }\alpha >0,
\label{2-0}
\end{equation}%
and the output corresponding to the controlled system has the property 
\begin{equation}
\left\Vert z\right\Vert _{L^{2}(\mathbb{R}_{+};Z)}=\left\Vert \Phi
(y,u)\right\Vert _{L^{2}(\mathbb{R}_{+};Z)}\leq \gamma \left\Vert
w\right\Vert _{L^{2}(\mathbb{R}_{+};W)}  \label{3}
\end{equation}%
for $\gamma $ a prescribed positive value representing the robust
performance index. We recall by a Datko's result (see \cite{Datko-72}) that (%
\ref{2-0}) is equivalent with $e^{(-A+B_{2}F)t}y\in L^{2}(\mathbb {R}_{+};H)\mathit{.}$
Results in the literature show that the feedback operator is constructed by
the relation 
\begin{equation}
F=-B_{2}^{\ast }P,  \label{4}
\end{equation}%
where the linear operator $P$ is the solution to the associated Riccati
equation 
\begin{equation}
A^{\ast }P+PA+P(B_{2}B_{2}^{\ast }-\gamma ^{-2}B_{1}B_{1}^{\ast
})P=C_{1}^{\ast }C_{1}.  \label{5}
\end{equation}%
Moreover, $P=P^{\ast }$ and $P\geq 0.$ This result follows under some
supplementary hypotheses and the operator $P$ has also some other properties
which we do not recall here since they will be not necessarily involved in
the next proofs. For more details on the $H^{\infty }$-optimal control
problem for parabolic systems we can refer the reader to \cite{GM-23-Hinf}.

There has been a large amount of literature on $H^{\infty }$-optimal control
problems for finite-dimensional systems, since 1980 (e.g., \cite{Basar}, 
\cite{Ichikawa-91}, \cite{Zhou-Doyle-Glover}), followed by papers addressing
infinite-dimensional systems with bounded control operators (see e.g., \cite%
{Bensoussan}, \cite{Ichikawa-96}), or the monograph \cite{Keulen} addressing
the $H^{\infty }$-control problem for a class of systems with unbounded
control operators. Papers \cite{VB-hyp-92} and \cite{VB-hyp-95} focused on $%
H^{\infty }$-boundary control under state feedback for hyperbolic-like
systems.

While in the finite-dimensional case the feedback operator can be found by
solving a matrix algebraic Riccati equation, its finding in the case of
infinite-dimensional systems requires to solve an infinite-dimensional
operator algebraic Riccati equation. For a general theory regarding the
Riccati equation in the infinite-dimensional case and its application we
refer the reader to the monographs \cite{IL-RT-1991} and \cite{IL-RT-1999}.
An operatorial method for solving the Riccati equation for $A^{\ast
}P+PA+P^{2}=F$ was given in \cite{VB-76}, Chapter II, 3.3, see also \cite%
{Brezis-1970}. A study for a very specific heat equation was done in \cite%
{BMR-2020}. An approximation theory of solutions to operator Riccati
equations for $H^{\infty }$- optimal control have been developed e.g., in 
\cite{Ito} and \cite{Morris} for bounded control operators.

In \cite{GM-23-Hinf} the $H^{\infty }$-optimal control abstract problem has
been solved for linear infinite dimensional systems of parabolic type,
relying on appropriate assumptions for parabolic operators and applications
have been done for equations with singular Hardy potentials. The proof was
essentially based on a variational technique, starting by solving the
differential game $\sup_{w\in L^{2}(\mathbb{R}_{+},W)}\inf_{u\in L^{2}(%
\mathbb{R}_{+},U)}\left\{ \frac{1}{2}\int_{0}^{\infty }\left( \left\Vert
z(t)\right\Vert _{Z}^{2}-\gamma \left\Vert w(t)\right\Vert _{W}^{2}\right)
dt\right\} .$

In what concerns the determination of the feedback operator by numerical
techniques, a sparse grid discretization for the numerical solution of the
algebraic Riccati equation was approached in \cite{Har-Kalm} by formulating
it as a nonlinear operator equation and by deriving the matrix equation for
the sparse grid discretization. The case of unbounded operators was treated
in \cite{Guo-Tan} under Dirichlet boundary control by employing the Galerkin
approximation, which generates a sequence of solvable finite-dimensional
systems that approximate the original infinite-dimensional system and by
demonstrating that the sequence of solutions to the associated
finite-dimensional algebraic Riccati equations converges in norm to the
solution of the infinite-dimensional operator algebraic Riccati equation.

This paper deals essentially with a direct approach of the existence of the
solution to the Riccati equation, completely different from that followed in 
\cite{GM-23-Hinf}, namely by using an operatorial method set up in the space
of Hilbert-Schmidt operators. We shall investigate the existence of a
solution to equation (\ref{5}) by writing it as 
\begin{equation}
A^{\ast }Py+PAy+P\Gamma Py=C_{1}^{\ast }C_{1}y  \label{6}
\end{equation}%
where $\Gamma :H\rightarrow H$ is defined by 
\begin{equation}
\Gamma :=B_{2}B_{2}^{\ast }-\gamma ^{-2}B_{1}B_{1}^{\ast }.  \label{7}
\end{equation}%
The outline of the paper is the following. After some preliminaries
introducing the appropriate functional framework for the study of the
existence of the solution to the Riccati equation, a first result is given
for in the case when the operator $A$ is coercive, in Theorem \ref{T2.1}. A
more challenging problem refers to $A\geq 0$ and this case is treated in
Theorem \ref{T3.1}. Relying on the existence of the solution to the Riccati
equation treated in the previous sections, Theorem \ref{T4.1} solves the $%
H^{\infty }$-optimal control problem, that is it shows that the control of
the form $F=-B_{2}^{\ast }P$ satisfies relation (\ref{3}). An example,
involving a parabolic equation with a singular term of Hardy type is given
for checking the application of the theoretical results regarding the
existence of the associated Riccati equation and this ends the paper. With
respect to the variational method followed in paper \cite{GM-23-Hinf}, the
approach considered here has an advantage because the proofs are
constructive, possibly indicating a method for calculating a numerical
solution.

\subsection{Preliminaries}

Let $V$ and $H$ be separable Hilbert spaces, with $V^{\prime }$ the dual of $%
V$, and set the variational triplet $V\subset H\subset V^{\prime }$ with
dense and continuous injections, assuming in addition that the embedding $V$
in $H$ is compact. We denote by a subscript $X$ the scalar product $(\cdot
,\cdot )_{X}$ and the norm $\left\Vert \cdot \right\Vert _{X}$ in a Hilbert
space $X$ and by $\left\langle \cdot ,\cdot \right\rangle _{V^{\prime },V}$
the pairing between $V^{\prime }$ and $V.$ We have 
\begin{equation}
\left\Vert y\right\Vert _{H}\leq \kappa _{1}\left\Vert y\right\Vert _{V},%
\mbox{ }\left\Vert y\right\Vert _{V^{\prime }}\leq \kappa _{2}\left\Vert
y\right\Vert _{H},\mbox{ }\kappa _{1},\mbox{ }\kappa _{2}>0.  \label{7-1}
\end{equation}

We denote by $L(H_{1},H_{2})$ the space of linear and continuous operators
from $H_{1}$ to $H_{2},$ $H_{i}$ being Hilbert spaces, $i=1,2$.

Let $J$ be the canonical isomorphism from $V$ to $V^{\prime },$ that is $%
\left\langle Ju,v\right\rangle _{V^{\prime },V}=(u,v)_{V}$ for $u,$ $v\in V,$
implying that $\left\langle Ju,u\right\rangle _{V^{\prime },V}=\left\Vert
u\right\Vert _{V}^{2}=\left\Vert Ju\right\Vert _{V^{\prime }}^{2}.$ Since $V$
is compact in $H$ then $J^{-1}$ belongs to $L(V^{\prime },V)\cap L(H,V)$ and
it is self-adjoint and compact on $H.$

In this section we consider $\{e_{n}\}_{n\geq 1}$ the orthonormal basis in
the Hilbert space $H,$ set such that 
\begin{equation}
J^{-1}e_{n}=\rho _{n}^{-2}e_{n},  \label{7-2}
\end{equation}%
that is $\{\rho _{n}^{2}\}_{n\geq 1}$ is the positive sequence of
eigenvalues for $J,$ increasingly ordered, 
\[
0<\rho ^{2} _{1}<\rho ^{2} _{2}<\mbox{ }... 
\]%
Such a basis exists due to the assumption that $V$ is compact in $H$. We
have $Je_{n}=\rho _{n}^{2}e_{n}$, and by (\ref{7-2}) it follows that $%
e_{n}\in V,$ too.

Then, $\widetilde{e_{n}}=\rho _{n}^{-1}e_{n}$ is an orthonormal basis in $V$
and $\widetilde{\widetilde{e_{n}}}=\rho _{n}e_{n}$ is an orthonormal basis
in $V^{\prime }.$ Indeed, we have 
\[
\left\Vert Je_{n}\right\Vert _{V^{\prime }}^{2}=\left\langle
Je_{n},e_{n}\right\rangle _{V^{\prime },V}=\rho _{n}^{2}\left\langle
e_{n},e_{n}\right\rangle _{V^{\prime },V}=\rho _{n}^{2}\left\Vert
e_{n}\right\Vert _{H}^{2}=\rho _{n}^{2},\mbox{ } 
\]%
hence%
\begin{equation}
\left\Vert J(\rho _{n}^{-1}e_{n})\right\Vert _{V^{\prime }}=\left\Vert \rho
_{n}^{-1}e_{n}\right\Vert _{V}=1  \label{7-3}
\end{equation}%
and%
\begin{equation}
\left\Vert e_{n}\right\Vert _{V^{\prime }}=\left\Vert J^{-1}e_{n}\right\Vert
_{V}=\left\Vert \rho _{n}^{-2}e_{n}\right\Vert _{V}=\rho _{n}^{-1}\left\Vert
\rho _{n}^{-1}e_{n}\right\Vert _{V}=\rho _{n}^{-1}.  \label{7-4}
\end{equation}

For a later use we give here the definition of a maximal monotone operator
on a Hilbert space $H$. Let $A:D(A)\subset H\rightarrow H$ be a nonlinear
operator on $H$ with the domain $D(A)$.

\begin{definition}
The operator $A$ is called monotone, if $(Ax_{1}-Ax_{2},x_{1}-x_{2})_{H}\geq
0,$ for any $x_{1},x_{2}\in D(A).$

The operator $A$ is called maximal monotone ($m$-accretive), if, in addition
to monotonicity, the operator satisfies the relation $R(I+A)=H,$ where $R$
is the range of the operator and $I$ is the identity operator on $H.$
\end{definition}

Now, we shall set the appropriate functional framework to prove existence
for the solution to the Riccati equation (\ref{6}).

\begin{definition}
Let $H_{i},$ $i=1,2,3$ be separable Hilbert spaces. The operator $%
P:H_{1}\rightarrow H_{2}$ is said to be Hilbert-Schmidt if for some complete
orthonormal basis $\{e_{n}\}_{n\geq 1}$ of $H_{1}$ we have 
\[
\sum\limits_{n=1}^{\infty }(Pe_{n},Pe_{n})_{H_{2}}<\infty . 
\]%
We denote by $(H_{1},H_{2})_{H.S.}$ the space of Hilbert-Schmidt operators
from $H_{1}$ to $H_{2}$ and recall that it is a Hilbert space with the
scalar product and norm 
\begin{equation}
(P,Q)_{H.S.}=\sum\limits_{n=1}^{\infty }(Pe_{n},Qe_{n})_{H_{2}},\mbox{ }%
\left\Vert P\right\Vert _{H.S.}^{2}=\sum\limits_{n=1}^{\infty }\left\Vert
Pe_{n}\right\Vert _{H_{2}}^{2},\mbox{ }P,\mbox{ }Q\in (H_{1},H_{2})_{H.S.}.
\label{11-0}
\end{equation}
\end{definition}

Next, we list some properties of the Hilbert-Schmidt operators (see e.g., 
\cite{D-S}, Chapter XI, see also \cite{VB-76}, Chapter II, 3.3):

(a) $\left\Vert P\right\Vert _{L(H_{1},H_{2})}\leq \left\Vert P\right\Vert
_{(H_{1},H_{2})_{H.S.}}$ for $P\in (H_{1},H_{2})_{H.S.},$

(b) If $P\in (H_{1},H_{2})_{H.S.}$ then $P$ is compact,

(c) If $P_{1}\in (H_{1},H_{2})_{H.S.}$ and $P_{2}\in (H_{2},H_{3})_{H.S.}$
then $P_{2}P_{1}\in (H_{1},H_{3})_{H.S.}$ and 
\[
\left\Vert P_{2}P_{1}\right\Vert _{(H_{1},H_{3})_{H.S.}}\leq \left\Vert
P_{2}\right\Vert _{(H_{2},H_{3})_{H.S.}}\left\Vert P_{1}\right\Vert
_{(H_{1},H_{2})_{H.S.}}, 
\]

(d) If $P_{1},$ $P_{2}\in (H,H)_{H.S.}$ and $P_{2}\geq 0$ then $%
(P_{2}P_{1},P_{1})_{H.S.}\geq 0$ and $(P_{1}P_{2},P_{1})_{H.S.}\geq 0.$

We still have 
\begin{equation}
(V^{\prime },H)_{H.S.}\subset (H,H)_{H.S.}\subset (V,H)_{H.S.},\mbox{ \ }%
(H,V)_{H.S.}\subset (H,H)_{H.S.}\subset (H,V^{\prime })_{H.S.}  \label{11-1}
\end{equation}%
and denote 
\begin{equation}
\mathcal{H}:=(H,H)_{H.S.},\mbox{ }\mathcal{V}:=(V^{\prime },H)_{H.S.}\cap
(H,V)_{H.S.}.  \label{11-2}
\end{equation}%
The dual space of $\mathcal{V},$ denoted $\mathcal{V}^{\prime }$ can be
identified with $(V,H)_{H.S.}+(H,V^{\prime })_{H.S.}$ and we have the
embeddings $\mathcal{V}\subset \mathcal{H}\subset \mathcal{V}^{\prime }$
with dense and continuous inclusions.

For the reader convenience we present the justification of the previous
inclusions (\ref{11-1}). For the rest of this section , let $\{e_{j}\}_{j}$
be the orthonormal basis in $H$ defined by (\ref{7-2}), implying that $%
e_{j}\in V$. We have 
\[
\left\Vert P\right\Vert _{(H,V^{\prime })_{H.S.}}^{2}=\sum_{j=1}^{\infty
}\left\Vert Pe_{j}\right\Vert _{V^{\prime }}^{2}\leq \kappa
_{2}^{2}\sum_{j=1}^{\infty }\left\Vert Pe_{j}\right\Vert _{H}^{2}=\kappa
_{2}^{2}\left\Vert P\right\Vert _{(H,H)_{H.S.}}^{2}, 
\]%
\[
\left\Vert P\right\Vert _{(H,H)_{H.S.}}^{2}=\sum_{j=1}^{\infty }\left\Vert
Pe_{j}\right\Vert _{H}^{2}\leq \kappa _{1}^{2}\sum_{j=1}^{\infty }\left\Vert
Pe_{j}\right\Vert _{V}^{2}=\kappa _{1}^{2}\left\Vert P\right\Vert
_{(H,V)_{H.S.}}^{2}, 
\]%
\[
\left\Vert P\right\Vert _{(V^{\prime },H)_{H.S.}}^{2}=\sum_{j=1}^{\infty
}\left\Vert P\widetilde{\widetilde{e_{j}}}\right\Vert
_{H}^{2}=\sum_{j=1}^{\infty }\rho _{j}^{2}\left\Vert Pe_{j}\right\Vert
_{H}^{2}\geq \rho _{1}^{2}\left\Vert P\right\Vert _{(H,H)_{H.S.}}^{2}, 
\]%
\[
\left\Vert P\right\Vert _{(V,H)_{H.S.}}^{2}=\sum_{j=1}^{\infty }\left\Vert P%
\widetilde{e_{j}}\right\Vert _{H}^{2}=\sum_{j=1}^{\infty }\rho
_{j}^{-2}\left\Vert Pe_{j}\right\Vert _{H}^{2}\leq \rho _{1}^{-2}\left\Vert
P\right\Vert _{(H,H)_{H.S.}}^{2}. 
\]%
Then,%
\begin{equation}
\left\Vert P\right\Vert _{\mathcal{V}^{\prime }}=\left\Vert P\right\Vert
_{(H,V^{\prime })_{H.S.}}+\left\Vert P\right\Vert _{(V,H)_{H.S.}}\leq
(\kappa _{2}+\rho _{1}^{-1})\left\Vert P\right\Vert _{\mathcal{H}},
\label{11-3}
\end{equation}%
\begin{equation}
\left\Vert P\right\Vert _{\mathcal{V}}=\left\Vert P\right\Vert _{(V^{\prime
},H)_{H.S.}}+\left\Vert P\right\Vert _{(H,V)_{H.S.}}\geq \left( \rho
_{1}+\kappa _{1}^{-1}\right) \left\Vert P\right\Vert _{\mathcal{H}}.
\label{11-4}
\end{equation}%
We define the pairings%
\begin{eqnarray*}
\left\langle P,Q\right\rangle _{(V,H)_{H.S.},(V^{\prime },H)_{H.S.}}
&=&\sum_{j=1}^{\infty }\left( P\widetilde{e_{j}},Q\widetilde{\widetilde{e_{j}%
}}\right) _{H},\mbox{ } \\
\left\langle P,Q\right\rangle _{(H,V^{\prime })_{H.S.},(H,V)_{H.S.}}
&=&\sum_{j=1}^{\infty }\left\langle Pe_{j},Qe_{j}\right\rangle _{V^{\prime
},V}.
\end{eqnarray*}

\begin{lemma}
\label{L1.2}Let $P\in \mathcal{V}$ be self-adjoint and $A:V\rightarrow
V^{\prime }$ be linear, continuous and coercive,%
\begin{equation}
\left\langle Ay,y\right\rangle _{V^{\prime },V}\geq \omega \left\Vert
y\right\Vert _{V}^{2},\mbox{ for all }y\in V,\mbox{ }\omega >0.  \label{8}
\end{equation}%
We have%
\begin{eqnarray}
\left\langle PA,P\right\rangle _{(V,H)_{H.S.},(V^{\prime },H)_{H.S.}} &\geq
&\omega \left\Vert P\right\Vert _{(V^{\prime },H)_{H.S.}}^{2}\mbox{ for all }%
P\in (V^{\prime },H)_{H.S.},  \label{20} \\
\left\langle A^{\ast }P,P\right\rangle _{(H,V^{\prime
})_{H.S.},(H,V)_{H.S.}} &\geq &\omega \left\Vert P\right\Vert
_{(H,V)_{H.S.}}^{2}\mbox{ for all }P\in (H,V)_{H.S.}.  \label{20-0}
\end{eqnarray}
\end{lemma}

\noindent \textbf{Proof.} If $y\in H,$ it is written $y=\sum_{j=1}^{\infty
}(y,e_{j})_{H}e_{j}$ and $Py=\sum_{j=1}^{\infty }(y,e_{j})_{H}Pe_{j}$. We
denote $a_{jk}:=\left\langle Ae_{j},e_{k}\right\rangle _{V^{\prime },V},$ $%
p_{jk}:=(Pe_{j},e_{k})_{H}$ and calculate%
\[
\left\langle A^{\ast }P,P\right\rangle _{(H,V^{\prime
})_{H.S.},(H,V)_{H.S.}}=\sum_{j=1}^{\infty }\left\langle A^{\ast
}Pe_{j},Pe_{j}\right\rangle _{V^{\prime },V}\geq \omega \sum_{j=1}^{\infty
}\left\Vert Pe_{j}\right\Vert _{V}^{2}=\omega \left\Vert P\right\Vert
_{(H,V)_{H.S.}}^{2}.\mbox{ } 
\]%
Then,%
\begin{eqnarray*}
&&\left\langle PA,P\right\rangle _{(V,H)_{H.S.},(V^{\prime
},H)_{H.S.}}=\sum_{j=1}^{\infty }\left( PA\widetilde{e_{j}},P\widetilde{%
\widetilde{e_{j}}}\right) _{H}=\sum_{j=1}^{\infty }(PAe_{j},Pe_{j})_{H} \\
&=&\sum_{j=1}^{\infty }\left( \sum_{k=1}^{\infty }\left\langle
Ae_{j},e_{k}\right\rangle _{V^{\prime },V}Pe_{k},\sum_{l=1}^{\infty
}(Pe_{j},e_{l})_{H}e_{l}\right) _{H} \\
&=&\sum_{j=1}^{\infty }\left( \sum_{k=1}^{\infty }\left\langle
Ae_{j},e_{k}\right\rangle _{V^{\prime },V}\sum_{l=1}^{\infty
}p_{jl}p_{kl}\right) =\sum_{l=1}^{\infty }\left\langle A\sum_{j=1}^{\infty
}e_{j}p_{jl},\sum_{k=1}^{\infty }p_{kl}e_{k}\right\rangle _{V^{\prime },V} \\
&=&\sum_{l=1}^{\infty }\left\langle A\sum_{j=1}^{\infty
}e_{j}p_{jl},\sum_{j=1}^{\infty }p_{jl}e_{j}\right\rangle _{V^{\prime
},V}\geq \omega \sum_{l=1}^{\infty }\left\Vert \sum_{j=1}^{\infty
}e_{j}p_{jl}\right\Vert _{V}^{2}=\omega \sum_{l=1}^{\infty }\left(
\sum_{j=1}^{\infty }e_{j}p_{jl},\sum_{i=1}^{\infty }e_{i}p_{il}\right) _{V}
\\
&=&\omega \sum_{l=1}^{\infty }\sum_{j=1}^{\infty }\sum_{i=1}^{\infty
}p_{jl}p_{il}(e_{j},e_{i})_{V}=\omega \sum_{l=1}^{\infty }\sum_{j=1}^{\infty
}p_{jl}^{2}\rho _{j}^{2}.
\end{eqnarray*}%
On the other hand,%
\begin{eqnarray*}
\left\Vert P\right\Vert _{(V^{\prime },H)_{H.S.}}^{2} &=&\sum_{j=1}^{\infty
}\rho _{j}^{2}(Pe_{j},Pe_{j})_{H}=\sum_{j=1}^{\infty }\rho _{j}^{2}\left(
\sum_{k=1}^{\infty }(Pe_{j},e_{k})_{H}e_{k},\sum_{k=1}^{\infty
}(Pe_{j},e_{k})_{H}e_{k}\right) _{H} \\
&=&\sum_{j=1}^{\infty }\sum_{k=1}^{\infty }\rho _{j}^{2}p_{jk}^{2},
\end{eqnarray*}%
which compared with the previous result leads to the first inequality in (%
\ref{20}).\hfill $\square $

\section{The coercive case}

\setcounter{equation}{0}

Let

\begin{equation}
B_{1}\in L(W,H),\mbox{ }B_{2}\in L(U,H),\mbox{ }C_{1}\in L(H,Z),\mbox{ }%
D_{1}\in L(U,Z)  \label{12-1}
\end{equation}%
and assume in this section that $A:V\rightarrow V^{\prime }$ is linear,
continuous and satisfies the coercivity condition (\ref{8}),%
\[
\left\langle Ay,y\right\rangle _{V^{\prime },V}\geq \omega \left\Vert
y\right\Vert _{V}^{2},\mbox{ for all }y\in V,\mbox{ }\omega >0. 
\]
Recall (\ref{11-2}).

\begin{theorem}
\label{T2.1}Let 
\begin{equation}
F\in \mathcal{H},\mbox{ }F=F^{\ast },\mbox{ }F\geq 0,  \label{12}
\end{equation}%
\begin{equation}
\Gamma \in \mathcal{H},\mbox{ }\Gamma =\Gamma ^{\ast },\mbox{ }\Gamma \geq 0.
\label{13}
\end{equation}%
Then, there exists a unique $P\in \mathcal{V}$ solution to the equation 
\begin{equation}
A^{\ast }P+PA+P\Gamma P=F,  \label{14}
\end{equation}%
having the properties 
\begin{equation}
P=P^{\ast },\mbox{ }P\geq 0.  \label{15}
\end{equation}
\end{theorem}

\noindent \textbf{Proof.} We define 
\begin{equation}
\mathcal{L}:D(\mathcal{L})\subset \mathcal{H}\rightarrow \mathcal{H},\mbox{ }%
D(\mathcal{L})=\{P\in \mathcal{V};\mbox{ }\mathcal{L}(P)\in \mathcal{H}\},%
\mbox{ }\mathcal{L}(P):=A^{\ast }P+PA,  \label{17}
\end{equation}%
\begin{equation}
D=\{P\in \mathcal{H};\mbox{ }P=P^{\ast },\mbox{ }P\geq 0\},  \label{16}
\end{equation}%
\begin{equation}
\mathcal{B}:D(\mathcal{B})=D\subset \mathcal{H}\rightarrow \mathcal{H},%
\mbox{
}\mathcal{B}(P):=P\Gamma P  \label{18}
\end{equation}%
and write equation (\ref{14}) as 
\begin{equation}
\mathcal{L}(P)+\mathcal{B}(P)=F.  \label{16-0}
\end{equation}%
The proof will be done in three steps. In the first two steps we check the
properties of the operators $\mathcal{L}$ and $\mathcal{B},$ and show that $%
\mathcal{L}$ is maximal monotone and that $\mathcal{B}$ is continuous and
monotone in $\mathcal{H}\times \mathcal{H}.$ In the last step the existence
of the solution to (\ref{16-0}) is proved by passing to the limit in an
approximating equation obtained by approximating the operator $\mathcal{B}.$

\noindent \textbf{Step 1.} We shall prove that the linear operator $\mathcal{%
L}$ is maximal monotone on $\mathcal{H}$ and the operator $(\mathcal{I}%
+\lambda \mathcal{L})^{-1}$ leaves invariant the domain $D.$ Here, $\mathcal{%
I}$ is the identity operator in $\mathcal{H},$ that is $\mathcal{I}:\mathcal{%
H}\rightarrow \mathcal{H},$ $\mathcal{I}(P)=P.$ Using the definition of the
space $\mathcal{V}$ and (\ref{20})-(\ref{20-0}) it follows immediately that 
\begin{equation}
\left\langle \mathcal{L(}P),P\right\rangle _{\mathcal{V}^{\prime }\mbox{,}%
\mathcal{V}}\geq \omega \left\Vert P\right\Vert _{\mathcal{V}}^{2}.
\label{21}
\end{equation}%
Then, we shall prove the invariance property 
\begin{equation}
(\mathcal{I}+\lambda \mathcal{L)}^{-1}D\subset D\mbox{ for all }\lambda >0.
\label{22}
\end{equation}%
Taking into account that the linear operator $P\rightarrow A^{\ast }P+PA$ is
continuous, monotone and coercive from $\mathcal{V}$ to $\mathcal{V}^{\prime
}$, it follows that it is surjective and so, $\mathcal{L},$ its restriction
to $\mathcal{H}$ is maximal monotone (see \cite{VB-2010}, p. 36, Corollary
2.2). Thus, for every $F\in \mathcal{H}$ the equation 
\begin{equation}
P+\lambda A^{\ast }P+\lambda PA=F  \label{23}
\end{equation}%
has a unique solution $P\in \mathcal{H}$. By multiplying (\ref{23}) scalarly
by $P$ we get $\lambda \omega \left\Vert P\right\Vert _{\mathcal{V}}^{2}\leq
\left\Vert F\right\Vert _{\mathcal{H}}^{2},$ hence $P\in \mathcal{V}$ and
still by (\ref{23}) we deduce that $\mathcal{L(}P)\in \mathcal{H},$ so that $%
P\in D(\mathcal{L}).$ Because the solution is unique it follows that $%
P=P^{\ast },$ since $P^{\ast }$ verifies (\ref{23}).

Next, we show that $P\in D$ if $F\in D.$ Since $P$ is compact (see property\
(b) of Hilbert-Schmidt operators) and self-adjoint it follows that it has a
sequence of eigenvectors $\{e_{j}^{P}\}_{j}$ which generates an orthonormal
basis in $H.$ The corresponding eigenvalues are denoted by $\lambda _{j}^{P}$%
, and we have 
\begin{equation}
Pe_{j}^{P}=\lambda _{j}^{P}e_{j}^{P},\mbox{ }y=\sum_{j=1}^{\infty
}y_{j}e_{j}^{P},\mbox{ }y_{j}=(y,e_{j}^{P})_{H}\mbox{ for all }y\in H.
\label{19}
\end{equation}%
Since $P\in \mathcal{V},$ it follows that $Pe_{j}^{P}\in V$ and so $%
e_{j}^{P}\in V.$ Equation (\ref{23}) can be equivalently written%
\[
Py+\lambda (A^{\ast }P+PA)y=Fy,\mbox{ for }y\in H. 
\]%
Taking $y=e_{j}^{P}$ and multiplying scalarly by $e_{j}^{P}$ we get 
\[
(Pe_{j}^{P},e_{j}^{P})_{H}+2\lambda \left\langle
Ae_{j}^{P},Pe_{j}^{P}\right\rangle _{V^{\prime
},V}=(Fe_{j}^{P},e_{j}^{P})_{H},\mbox{ } 
\]%
Then, using (\ref{19}) we obtain 
\[
(\lambda _{j}^{P}e_{j}^{P},e_{j}^{P})_{H}+2\lambda \lambda
_{j}^{P}\left\langle Ae_{j}^{P},e_{j}^{P}\right\rangle _{V^{\prime
},V}=(Fe_{j}^{P},e_{j}^{P})_{H} 
\]%
and this yields 
\begin{equation}
\lambda _{j}^{P}=\left( 1+2\lambda \left\langle
Ae_{j}^{P},e_{j}^{P}\right\rangle _{V^{\prime },V}\right)
^{-1}(Fe_{j}^{P},e_{j}^{P})_{H},  \label{23-0}
\end{equation}%
which is nonnegative because $\left\langle Ae_{j}^{P},e_{j}^{P}\right\rangle
_{V^{\prime },V}\geq 0$ and $F\geq 0.$ In conclusion, $\lambda _{j}^{P}\geq
0 $ and so $P\geq 0.$

\noindent \textbf{Step 2.} We have 
\begin{equation}
\left\Vert \mathcal{B}(P)\right\Vert _{\mathcal{H}}=\left\Vert P\Gamma
P\right\Vert _{\mathcal{H}}\leq \left\Vert \Gamma \right\Vert _{\mathcal{H}%
}\left\Vert P\right\Vert _{\mathcal{H}}^{2}\mbox{ for }P\in \mathcal{H},
\label{24}
\end{equation}%
which follows by the property (c), and show then that $\mathcal{B}$ is
monotone in $\mathcal{H}\times \mathcal{H}.$ For $P_{1},$ $P_{2}\in D,$ we
have to prove that 
\begin{eqnarray}
E &:&=(\mathcal{B}(P_{1})-\mathcal{B}(P_{2}),P_{1}-P_{2})_{\mathcal{H}%
}=(P_{1}\Gamma P_{1}-P_{2}\Gamma P_{2},P_{1}-P_{2})_{\mathcal{H}}  \label{25}
\\
&=&((P_{1}-P_{2})\Gamma P_{1},P_{1}-P_{2})_{\mathcal{H}}+(P_{2}\Gamma
(P_{1}-P_{2}),P_{1}-P_{2})_{\mathcal{H}}\geq 0.  \nonumber
\end{eqnarray}%
It is obvious that if the sufficient conditions 
\begin{equation}
P\Gamma \geq 0\mbox{ and }\Gamma P\geq 0\mbox{ }  \label{26}
\end{equation}%
hold for all $P\in D$, then $E\geq 0,$ due to hypothesis (d).

We have to prove that, under the hypothesis (\ref{13}), $\Gamma \geq 0,$
conditions (\ref{26}) hold true. First we show that 
\begin{equation}
(\Gamma Py,y)_{H}\geq 0\mbox{ for all }y\in H.  \label{27}
\end{equation}%
As before, let $\{e_{j}^{P}\}_{j}$ be a sequence of eigenvectors of $P$
which generates an orthonormal basis in $H.$ The corresponding eigenvalues
are denoted by $\lambda _{j}^{P}$, hence $Pe_{j}^{P}=\lambda
_{j}^{P}e_{j}^{P}$, $j=1,2...$ and we note that 
\begin{equation}
\lambda _{j}^{P}\geq 0\mbox{ for all }j=1,2,...  \label{28}
\end{equation}%
We write again $y$ and $Py$ as in (\ref{19}), and replace in (\ref{27}). We
have 
\begin{eqnarray}
I_{1} &:&=(Py,\Gamma y)_{H}=\left( \sum_{j=1}^{\infty
}y_{j}Pe_{j}^{P},\sum_{k=1}^{\infty }y_{k}\Gamma e_{k}^{P}\right)
_{H}=\sum_{j,k=1}^{\infty }y_{j}y_{k}(Pe_{j}^{P},\Gamma e_{k}^{P})_{H}
\label{30} \\
&=&\sum_{j,k=1}^{\infty }y_{j}y_{k}\lambda _{j}^{P}(\Gamma
e_{k}^{P},e_{j}^{P})_{H}.  \nonumber
\end{eqnarray}%
By hypothesis (\ref{13}), the operator $\Gamma $ is positive semidefinite,
so that the associate quadratic form 
\[
I_{\Gamma }=\sum_{j,k=1}^{\infty }y_{j}y_{k}(\Gamma e_{k}^{P},e_{j}^{P})_{H} 
\]%
is positive semidefinite. This is equivalent with the fact that its infinite
matrix representation $M_{\Gamma }=\left\Vert (\Gamma
e_{k}^{P},e_{j}^{P})_{H}\right\Vert _{j,k\geq 1}$ is positive semidefinite,
and this takes place if all principal minors 
\begin{equation}
\det \left\Vert 
\begin{tabular}{lllllll}
$(\Gamma e_{1},e_{1})$ & ... & $(\Gamma e_{1},e_{j-1})$ & $(\Gamma
e_{1},e_{j+1})$ & ... & $(\Gamma e_{1},e_{n})$ & ... \\ 
... &  &  &  &  &  &  \\ 
$(\Gamma e_{j-1},e_{1})$ & ... & $(\Gamma e_{j-1},e_{j-1})$ & ... &  &  & 
\\ 
$(\Gamma e_{j+1},e_{1})$ & ... & $(\Gamma e_{j+1},e_{j-1})$ & ... &  &  & 
\\ 
$(\Gamma e_{n},e_{1})$ & ... &  &  &  & $(\Gamma e_{n},e_{n})$ &  \\ 
... &  &  &  &  & ... & ...%
\end{tabular}%
\right\Vert  \label{31-1}
\end{equation}%
are nonnegative.

Now, it is nothing else to do than noticing that by multiplying each row of
the principal minors of the infinite matrix $M_{\Gamma }$ by $\lambda
_{j}^{P}\geq 0,$ it follows that they remain nonnegative. But these are
exactly the principal minors of the infinite matrix $M_{1}:=\left\Vert
\lambda _{j}^{P}(\Gamma e_{k}^{P},e_{j}^{P})_{H}\right\Vert _{j,k\geq 1}$
and this implies that it is positive semidefinite and then the associated
quadratic form $I_{1}\geq 0.$

Similarly we show that $P\Gamma \geq 0$ for all $P\in D.$ In conclusion (\ref%
{26}) holds, which implies that $E\geq 0$ and so we have ended the proof of
the monotonicity of the operator $\mathcal{B}.$

The proof of Theorem \ref{T2.1} is continued by showing still that 
\begin{equation}
(\mathcal{I}+\lambda \mathcal{B)}^{-1}D\subset D.  \label{32}
\end{equation}%
To this end we have to show that the equation 
\begin{equation}
P+\lambda \mathcal{B}(P)=F,\mbox{ for any }F\in D  \label{33}
\end{equation}%
has a solution $P\in D.$ By (\ref{12}) it follows that $F$ is compact too
and so it has a sequence of eigenvectors $\{f_{j}\}_{j}$ which generates an
orthonormal basis in $H.$ We denote the eigenvalues of $F$ by $\lambda
_{j}^{F}$ and so we have 
\[
Ff_{j}=\gamma _{j}^{F}f_{j},\mbox{ }j\geq 1,\mbox{ with }\gamma _{j}^{F}\geq
0 
\]%
due to the fact that $F\geq 0.$

For the beginning, let us consider the equation 
\begin{equation}
\lambda g_{jj}\gamma _{j}^{2}+\gamma _{j}-\gamma _{j}^{F}=0,\mbox{ for all }%
j\geq 1,  \label{40}
\end{equation}%
where $g_{jj}:=(\Gamma e_{j},e_{j})_{H}.$ This equation has a nonnegative
solution for all $j\geq 1,$ 
\begin{equation}
\gamma _{j}=\left\{ 
\begin{tabular}{ll}
$\frac{-1+\sqrt{1+4\lambda g_{jj}\gamma _{j}^{F}}}{2\lambda g_{jj}}$ & $%
\mbox{if }g_{jj}\neq 0,\mbox{ }$ \\ 
$\gamma _{j}^{F},$ & if $g_{jj}=0$%
\end{tabular}%
\right.  \label{41}
\end{equation}%
which has the properties 
\begin{equation}
0\leq \gamma _{j}\leq \frac{2\gamma _{j}^{F}}{1+\sqrt{1+4\lambda
g_{jj}\gamma _{j}^{F}}}\leq \gamma _{j}^{F}.  \label{41-0}
\end{equation}%
Now, we write 
\begin{equation}
y=\sum_{j=1}^{\infty }y_{j}f_{j},\mbox{ }y_{j}=(y,f_{j})_{H},\mbox{ for all }%
y\in H  \label{34}
\end{equation}%
and define 
\begin{equation}
Py:=\sum_{j=1}^{\infty }\gamma _{j}y_{j}f_{j},\mbox{ for all }y\in H
\label{35}
\end{equation}%
which implies that $\gamma _{j}$ are the eigenvalues of $P$ in the expansion
with respect to $\{f_{j}\}_{j},$ 
\begin{equation}
Pf_{j}=\gamma _{j}f_{j},\mbox{ }j\geq 1.  \label{36}
\end{equation}%
We prove that $P$ given by (\ref{35}) satisfies (\ref{33}). First we
calculate 
\begin{equation}
\mathcal{B}(P)y=P\Gamma Py=P\Gamma \left( \sum_{j=1}^{\infty }\gamma
_{j}y_{j}f_{j}\right) =P\left( \sum_{j=1}^{\infty }\gamma _{j}y_{j}\Gamma
f_{j}\right) =\sum_{k=1}^{\infty }\left( \sum_{j=1}^{\infty }\gamma
_{j}y_{j}\Gamma f_{j},f_{k}\right) _{H}\gamma _{k}f_{k}  \label{36-0}
\end{equation}%
and replace in (\ref{33}) calculated in $y$. We have 
\begin{equation}
\sum_{j=1}^{\infty }\gamma _{j}y_{j}f_{j}+\lambda \sum_{j,k=1}^{\infty
}\gamma _{j}\gamma _{k}y_{j}g_{jk}f_{k}=\sum_{j=1}^{\infty }y_{j}\gamma
_{j}^{F}f_{j},\mbox{ for all }y=\{y_{j}\}_{\dot{j}\geq 1},  \label{37}
\end{equation}%
where 
\begin{equation}
g_{jk}:=(\Gamma f_{j},f_{k})_{H}.  \label{38}
\end{equation}%
Since (\ref{37}) takes place for all $y=\{y_{j}\}_{\dot{j}\geq 1}$ we deduce
the equation%
\begin{equation}
\gamma _{j}f_{j}+\lambda \sum_{k=1}^{\infty }\gamma _{j}\gamma
_{k}g_{jk}f_{k}=\gamma _{j}^{F}f_{j},\mbox{ for all }j\geq 1,  \label{39}
\end{equation}%
whence by multiplying by $f_{j}$ we finally obtain indeed that $\gamma _{j}$
is the solution to (\ref{40}) for each $j.$ Then, 
\[
\left\Vert P\right\Vert _{\mathcal{H}}^{2}=\sum\limits_{j=1}^{\infty
}\left\Vert Pf_{j}\right\Vert _{H}^{2}=\sum\limits_{j=1}^{\infty }|\gamma
_{j}|^{2}\left\Vert f_{j}\right\Vert _{H}^{2}\leq \sum\limits_{j=1}^{\infty
}|\gamma _{j}|^{2}\leq \sum\limits_{j=1}^{\infty }|\gamma _{j}^{F}|^{2}\leq
\sum\limits_{j=1}^{\infty }\left\Vert Ff_{j}\right\Vert _{H}^{2}\left\Vert
f_{j}\right\Vert _{H}^{2}\leq \left\Vert F\right\Vert _{\mathcal{H}}^{2}. 
\]%
We conclude that $P$ the unique solution to (\ref{33}) belong to $\mathcal{H}%
,$ $P\geq 0$ and $P=P^{\ast }.$ Thus, $(\mathcal{I}+\lambda \mathcal{B)}%
^{-1} $ leaves invariant $D$ and%
\begin{equation}
\left\Vert (\mathcal{I}+\lambda \mathcal{B)}^{-1}F\right\Vert _{\mathcal{H}%
}\leq \left\Vert F\right\Vert _{\mathcal{H}},\mbox{ for }F\in D.
\label{41-1}
\end{equation}

\noindent \textbf{Step 3.} Now we shall prove the existence of the solution
to equation (\ref{16-0}). For $\lambda >0$ let 
\begin{equation}
\mathcal{B}_{\lambda }P:=\lambda ^{-1}(\mathcal{I}-(\mathcal{I}+\lambda 
\mathcal{B})^{-1})P=\mathcal{B}((\mathcal{I}+\lambda \mathcal{B})^{-1}P),%
\mbox{ for
all }P\in \mathcal{H}  \label{41-2}
\end{equation}%
define the Yosida approximation of $\mathcal{B}$ and consider the
approximating equation%
\begin{equation}
\mathcal{L}(P_{\lambda })+\mathcal{B}_{\lambda }(P_{\lambda })=F.  \label{43}
\end{equation}%
The equation can be equivalently rewritten 
\begin{equation}
P_{\lambda }=\lambda \left( \mathcal{I}+\lambda \mathcal{L}\right)
^{-1}F+\left( \mathcal{I}+\lambda \mathcal{L}\right) ^{-1}(\mathcal{I}%
+\lambda \mathcal{B)}^{-1}P_{\lambda }.  \label{44}
\end{equation}%
We consider the mapping $\mathcal{T}:D\rightarrow \mathcal{H},$%
\[
\mathcal{T}P_{\lambda }=\lambda \left( \mathcal{I}+\lambda \mathcal{L}%
\right) ^{-1}F+\left( \mathcal{I}+\lambda \mathcal{L}\right) ^{-1}(\mathcal{I%
}+\lambda \mathcal{B)}^{-1}P_{\lambda }, 
\]%
where $D$ is organized as a metric space with the distance $%
d(P_{1},P_{2})=\left\Vert P_{1}-P_{2}\right\Vert _{\mathcal{H}}$ for all $%
P_{1},$ $P_{2}\in \mathcal{H}.$ By (\ref{22}) and (\ref{32}) it follows that 
$\mathcal{T}$ maps $D$ into $D.$ It remains to show that $\mathcal{T}$ is a
contraction.

We recall (\ref{41-1}) and show that the operator $\left( \mathcal{I}%
+\lambda \mathcal{L}\right) ^{-1}$ is Lipschitz with the constant $%
1/(1+\lambda \omega \left( \rho _{1}+\kappa _{1}^{-1}\right) ^{2})<1.$
Indeed, by multiplying scalarly by $P$ the equation 
\[
P+\lambda \mathcal{L}(P)=Z\in \mathcal{H}, 
\]%
and using (\ref{11-4}) we have 
\[
\left( 1+\lambda \omega \left( \rho _{1}+\kappa _{1}^{-1}\right) ^{2}\right)
\left\Vert P\right\Vert _{\mathcal{H}}^{2}\leq \left\Vert P\right\Vert _{%
\mathcal{H}}^{2}+\lambda \omega \left\Vert P\right\Vert _{\mathcal{V}%
}^{2}\leq \left\Vert Z\right\Vert _{\mathcal{H}}\left\Vert P\right\Vert _{%
\mathcal{H}}, 
\]%
hence 
\[
\left\Vert \left( \mathcal{I}+\lambda \mathcal{L}\right) ^{-1}Z\right\Vert _{%
\mathcal{H}}\leq \frac{1}{1+\lambda \omega \left( \rho _{1}+\kappa
_{1}^{-1}\right) ^{2}}\left\Vert Z\right\Vert _{\mathcal{H}}. 
\]%
Therefore, we have 
\begin{eqnarray*}
\left\Vert \mathcal{T}P_{\lambda ,1}-\mathcal{T}P_{\lambda ,2}\right\Vert _{%
\mathcal{H}} &\leq &\frac{1}{1+\lambda \omega \left( \rho _{1}+\kappa
_{1}^{-1}\right) ^{2}}\left\Vert (\mathcal{I}+\lambda \mathcal{B)}%
^{-1}P_{\lambda ,1}-(\mathcal{I}+\lambda \mathcal{B)}^{-1}P_{\lambda
,2}\right\Vert _{\mathcal{H}} \\
&\leq &\frac{1}{1+\lambda \omega \left( \rho _{1}+\kappa _{1}^{-1}\right)
^{2}}\left\Vert P_{\lambda ,1}-P_{\lambda ,2}\right\Vert _{\mathcal{H}}
\end{eqnarray*}%
and so it follows that $\mathcal{T}$ has a fixed point and equation (\ref{44}%
) has a unique solution $P_{\lambda }\in \mathcal{H},$ $P_{\lambda
}=P_{\lambda }^{\ast }\geq 0.$

Then, by multiplying (\ref{44}) scalarly in $\mathcal{H}$ by $P_{\lambda },$
taking into account (\ref{21}) and the fact that $(\mathcal{B}_{\lambda
}(P_{\lambda }),P_{\lambda })_{\mathcal{H}}\geq 0,$ we get%
\[
\omega \left( \rho _{1}+\kappa _{1}^{-1}\right) ^{2}\left\Vert P_{\lambda
}\right\Vert _{\mathcal{H}}^{2}\leq \omega \left\Vert P_{\lambda
}\right\Vert _{\mathcal{V}}^{2}\leq \left\Vert F\right\Vert _{\mathcal{H}%
}\left\Vert P_{\lambda }\right\Vert _{\mathcal{H}}. 
\]%
We deduce that 
\begin{equation}
\left\Vert P_{\lambda }\right\Vert _{\mathcal{H}}\leq C\left\Vert
F\right\Vert _{\mathcal{H}},\mbox{ }\left\Vert P_{\lambda }\right\Vert _{%
\mathcal{V}}\leq C\left\Vert F\right\Vert _{\mathcal{H}}  \label{45}
\end{equation}%
where $C$ denotes several positive constants independent of $\lambda .$ By (%
\ref{41-1}) we still have 
\begin{equation}
\left\Vert (\mathcal{I}+\lambda \mathcal{B)}^{-1}P_{\lambda }\right\Vert _{%
\mathcal{H}}\leq \left\Vert P_{\lambda }\right\Vert _{\mathcal{H}}\leq C
\label{47}
\end{equation}%
and since $\mathcal{B}_{\lambda }(P_{\lambda })=\mathcal{B}(\mathcal{I}%
+\lambda \mathcal{B)}^{-1}P_{\lambda }$ and 
\[
\left\Vert \mathcal{B}(\mathcal{I}+\lambda \mathcal{B)}^{-1}P_{\lambda
}\right\Vert _{\mathcal{H}}\leq \left\Vert \Gamma \right\Vert _{\mathcal{H}%
}\left\Vert (\mathcal{I}+\lambda \mathcal{B)}^{-1}P_{\lambda }\right\Vert _{%
\mathcal{H}}^{2} 
\]%
by (\ref{24}), it follows that%
\begin{equation}
\left\Vert \mathcal{B}_{\lambda }(P_{\lambda })\right\Vert _{\mathcal{H}%
}\leq C.  \label{49}
\end{equation}%
Moreover, we have%
\begin{equation}
\left\Vert P_{\lambda }-(\mathcal{I}+\lambda \mathcal{B)}^{-1}P_{\lambda
}\right\Vert _{\mathcal{H}}\leq \lambda C  \label{48}
\end{equation}%
by (\ref{41-2}) and 
\[
\left\Vert \mathcal{L}(P_{\lambda })\right\Vert _{\mathcal{H}}\leq C, 
\]%
by equation (\ref{43}).

Next, we prove that $(P_{\lambda })_{\lambda }$ is Cauchy in $\mathcal{H}.$
We\ recall (\ref{41-2}) and multiply the equation 
\[
\mathcal{L}(P_{\lambda })-\mathcal{L}(P_{\mu })+\mathcal{B}(\mathcal{I}%
+\lambda \mathcal{B)}^{-1}P_{\lambda }-\mathcal{B}(\mathcal{I}+\lambda 
\mathcal{B)}^{-1}P_{\mu }=0 
\]%
scalarly by $P_{\lambda }-P_{\mu }.$ We write 
\begin{eqnarray*}
&&(\mathcal{L}(P_{\lambda })-\mathcal{L}(P_{\mu }),P_{\lambda }-P_{\mu })_{%
\mathcal{H}} \\
&&+(\mathcal{B}(\mathcal{I}+\lambda \mathcal{B)}^{-1}P_{\lambda }-\mathcal{B}%
(\mathcal{I}+\lambda \mathcal{B)}^{-1}P_{\mu },(\mathcal{I}+\lambda \mathcal{%
B)}^{-1}P_{\lambda }-(\mathcal{I}+\lambda \mathcal{B)}^{-1}P_{\mu })_{%
\mathcal{H}} \\
&=&-(\mathcal{B}(\mathcal{I}+\lambda \mathcal{B)}^{-1}P_{\lambda }-\mathcal{B%
}(\mathcal{I}+\lambda \mathcal{B)}^{-1}P_{\mu },P_{\lambda }-(\mathcal{I}%
+\lambda \mathcal{B)}^{-1}P_{\lambda }-(\mathcal{P}_{\mu }-(\mathcal{I}%
+\lambda \mathcal{B)}^{-1}P_{\mu }))_{\mathcal{H}}.
\end{eqnarray*}%
Recalling (\ref{21}) and that $E$ in (\ref{25}) is nonnegative, we obtain 
\begin{eqnarray*}
&&\omega \left( \rho _{1}+\kappa _{1}^{-1}\right) ^{2}\left\Vert P_{\lambda
}-P_{\mu }\right\Vert _{\mathcal{H}}^{2}\leq \left\Vert \mathcal{B}(\mathcal{%
I}+\lambda \mathcal{B)}^{-1}P_{\lambda }-\mathcal{B}(\mathcal{I}+\lambda 
\mathcal{B)}^{-1}P_{\mu }\right\Vert _{\mathcal{H}}\times \\
&&\left( \left\Vert P_{\lambda }-(\mathcal{I}+\lambda \mathcal{B)}%
^{-1}P_{\lambda }\right\Vert _{\mathcal{H}}+\left\Vert P_{\mu }-(\mathcal{I}%
+\lambda \mathcal{B)}^{-1}P_{\mu }\right\Vert _{\mathcal{H}}\right) \\
&\leq &(\lambda +\mu )C
\end{eqnarray*}%
showing that $(P_{\lambda })_{\lambda }$ is Cauchy. By (\ref{45})-(\ref{49})
and the latter property, we conclude that 
\[
P_{\lambda }\rightarrow P\mbox{ strongly in }\mathcal{H}\mbox{ }(%
\mbox{that
is convergence in norm)},\mbox{ and weakly in }\mathcal{V}. 
\]%
Then, (\ref{43}) and a similar calculus as in (\ref{25}) yield 
\begin{eqnarray*}
&&\left\Vert \mathcal{B}_{\lambda }(P_{\lambda })-\mathcal{B}_{\lambda
}(P)\right\Vert _{\mathcal{H}}=\left\Vert \mathcal{B}(\mathcal{I}+\lambda 
\mathcal{B)}^{-1}P_{\lambda }-\mathcal{B}(\mathcal{I}+\lambda \mathcal{B)}%
^{-1}P\right\Vert _{\mathcal{H}} \\
&\leq &\left\Vert \left( (\mathcal{I}+\lambda \mathcal{B)}^{-1}P_{\lambda }-(%
\mathcal{I}+\lambda \mathcal{B)}^{-1}P\right) \Gamma (\mathcal{I}+\lambda 
\mathcal{B)}^{-1}P_{\lambda }\right\Vert _{\mathcal{H}} \\
&&+\left\Vert \left( (\mathcal{I}+\lambda \mathcal{B)}^{-1}P\right) \Gamma
\left( (\mathcal{I}+\lambda \mathcal{B)}^{-1}P_{\lambda }-(\mathcal{I}%
+\lambda \mathcal{B)}^{-1}P\right) \right\Vert _{\mathcal{H}}.
\end{eqnarray*}%
Now, 
\begin{eqnarray*}
&&\left\Vert (\mathcal{I}+\lambda \mathcal{B)}^{-1}P_{\lambda }-(\mathcal{I}%
+\lambda \mathcal{B)}^{-1}P\right\Vert _{\mathcal{H}} \\
&\leq &\left\Vert (\mathcal{I}+\lambda \mathcal{B)}^{-1}P_{\lambda
}-P_{\lambda }\right\Vert _{\mathcal{H}}+\left\Vert P_{\lambda
}-P\right\Vert _{\mathcal{H}}+\left\Vert P-(\mathcal{I}+\lambda \mathcal{B)}%
^{-1}P\right\Vert _{\mathcal{H}}
\end{eqnarray*}%
which tends to $0$ as $\lambda \rightarrow 0.$ Next, 
\begin{eqnarray*}
&&\left\Vert \mathcal{B}_{\lambda }(P)-\mathcal{B}(P)\right\Vert _{\mathcal{H%
}}=\left\Vert \mathcal{B}\left( (I+\lambda \mathcal{B})^{-1}P\right) -%
\mathcal{B}(P)\right\Vert _{\mathcal{H}} \\
&=&\left\Vert ((\mathcal{I}+\lambda \mathcal{B})^{-1}P-P)\Gamma (\mathcal{I}%
+\lambda \mathcal{B})^{-1}P\right\Vert _{\mathcal{H}}+\left\Vert P\Gamma ((%
\mathcal{I}+\lambda \mathcal{B})^{-1}P-P)\right\Vert _{\mathcal{H}}\leq
\lambda C.
\end{eqnarray*}%
Then, writing 
\[
\left\Vert \mathcal{B}_{\lambda }(P_{\lambda })-\mathcal{B}(P)\right\Vert _{%
\mathcal{H}}=\left\Vert \mathcal{B}_{\lambda }(P_{\lambda })-\mathcal{B}%
_{\lambda }(P)\right\Vert _{\mathcal{H}}+\left\Vert \mathcal{B}_{\lambda
}(P)-\mathcal{B}(P)\right\Vert _{\mathcal{H}} 
\]%
we deduce that 
\[
\mathcal{B}_{\lambda }(P_{\lambda })\rightarrow \mathcal{B}(P)%
\mbox{
strongly in }\mathcal{H}, 
\]%
and then,%
\[
\mathcal{L}(P_{\lambda })\rightarrow \eta =F\mathcal{-B}(P)%
\mbox{ strongly
in }\mathcal{H}, 
\]%
by equation (\ref{43}). Since $\mathcal{L}$ is closed it follows that $\eta =%
\mathcal{L}(P)$, which means, in conclusion, that $P$ satisfies the equation 
$\mathcal{L}(P)+\mathcal{B}(P)=F.$

\noindent The solution is unique due to the strongly monotonicity of the
operator $\mathcal{L}+\mathcal{B}.$ This ends the proof.\hfill $\square $

\begin{corollary}
\label{C2.2}Under the hypotheses of Theorem \ref{T2.1}, let $C_{1}\in 
\mathcal{H}.$ Then, equation (\ref{6}) (and the corresponding equation (\ref%
{5})) has a unique solution $P,$ with $P\geq 0$ and self-adjoint.
\end{corollary}

It is nothing else to do than applying Theorem \ref{T2.1} with $%
F=C_{1}^{\ast }C_{1}.$ Then, $F\in $ $\mathcal{H}$ by property (c) in
Section 1.1 and $F\geq 0$ since $(Fy,y)_{H}=\left\Vert C_{1}y\right\Vert
_{H}^{2}\geq 0.$

\section{The noncoercive case}

\setcounter{equation}{0}

This section is dedicated to the noncoercive case, under the assumptions (%
\ref{12-1}). In this case, we shall strengthen the hypotheses for $\Gamma $
and $C_{1},$ imposing that $\Gamma $ is coercive and $C_{1}$ belongs to $%
\mathcal{V}$ and has the kernel consisting in the element $0.$

\begin{theorem}
\label{T3.1}Let%
\begin{equation}
A:V\rightarrow V^{\prime },\mbox{ }\left\langle Ay,y\right\rangle
_{V^{\prime },V}\geq 0,\mbox{ for all }y\in V,  \label{50}
\end{equation}%
\begin{equation}
C_{1}\in (V^{\prime },H)_{H.S.},\mbox{ }\ker C_{1}=\{0\},  \label{51}
\end{equation}%
\begin{equation}
\Gamma \in \mathcal{H},\mbox{ }\Gamma =\Gamma ^{\ast },\mbox{ }(\Gamma
y,y)_{H}\geq g_{0}\left\Vert y\right\Vert _{H}^{2},\mbox{ for any }y\in H,%
\mbox{ }g_{0}>0.  \label{52}
\end{equation}%
Then, there exists a unique $P\in \mathcal{V}$, positive and self-adjoint,
solution to the equation 
\begin{equation}
A^{\ast }P+PA+P\Gamma P=C_{1}^{\ast }C_{1}.  \label{53}
\end{equation}
\end{theorem}

\noindent \textbf{Proof.} Let us consider the canonical isomorphism $%
J:V\rightarrow V^{\prime },$ define the approximation $A_{\omega }$ of $A$
by 
\[
A_{\omega }:V\rightarrow V^{\prime },\mbox{ }A_{\omega }=A+\omega J%
\mbox{
for }\omega >0, 
\]%
and denote 
\[
\mathcal{L}_{\omega }(P)=A_{\omega }^{\ast }P+P_{\omega }A,\mbox{ }\mathcal{L%
}_{\omega }:D(\mathcal{L}_{\omega })=\{P\in \mathcal{V};\mbox{ }\mathcal{L}%
_{\omega }(P)\in H\}\subset \mathcal{H}\rightarrow \mathcal{H}. 
\]%
We recall the notation%
\[
\mathcal{L}(P)=A^{\ast }P+PA 
\]%
with $\mathcal{L}(P)$ defined in (\ref{17}), but corresponding to $A\geq 0$.
We consider equation (\ref{53}) for $A_{\omega },$ 
\begin{equation}
\mathcal{L}_{\omega }(P_{\omega })+P_{\omega }\Gamma P_{\omega }=F,\mbox{ }%
F:=C_{1}^{\ast }C_{1}.  \label{54}
\end{equation}%
By (\ref{51}) it follows that $F$ is positive, self-adjoint and $F\in 
\mathcal{H}.$

Since $A_{\omega }$ satisfies hypothesis (\ref{8}), $\left\langle A_{\omega
}y,y\right\rangle _{V^{\prime },V}\geq \omega \left\Vert y\right\Vert
_{V}^{2}$, it follows by Theorem \ref{T2.1} that equation (\ref{54}) has a
unique solution $P_{\omega }\in D(\mathcal{L}_{\omega }\mathcal{)\subset V}$%
, with $P_{\omega }\in D,$ that is $P_{\omega }=P_{\omega }^{\ast }$ and $%
P_{\omega }\geq 0$. Moreover, 
\begin{equation}
\left\Vert P_{\omega }\right\Vert _{\mathcal{H}}\leq \frac{1}{\sqrt{g_{0}}}%
\left\Vert C_{1}\right\Vert _{\mathcal{H}}.  \label{55}
\end{equation}%
We shall prove this assertion. We consider the orthonormal basis $%
\{e_{j}\}_{j}$ in $H,$ given by (\ref{7-2}), corresponding to the canonical
isomorphism $J:V\rightarrow V^{\prime },$ $Je_{j}=\rho _{j}^{2}e_{j},$ write 
\[
y=\sum_{j=1}^{\infty }y_{j}e_{j},\mbox{ for all }y\in H 
\]%
and define 
\[
P_{\omega }y=\sum_{j=1}^{\infty }\gamma _{j}^{\omega }y_{j}e_{j},%
\mbox{ for
all }y\in H. 
\]%
It follows that $P_{\omega }e_{j}=\gamma _{j}^{\omega }e_{j}.$ We prove that
indeed, there exist $\gamma _{j}^{\omega }\geq 0$ such that $P_{\omega }$
defined in this way verifies equation (\ref{54}). We calculate $\mathcal{B}%
(P)y$ similarly as by (\ref{36-0}) and write 
\[
(A_{\omega }^{\ast }P_{\omega }+P_{\omega }A_{\omega })y=\sum_{j=1}^{\infty
}y_{j}(A_{\omega }^{\ast }P_{\omega }+P_{\omega }A_{\omega })e_{j}. 
\]%
We apply both terms in (\ref{54}) to $e_{j}$ and deduce 
\[
(A_{\omega }^{\ast }P_{\omega }+P_{\omega }A_{\omega
})e_{j}+\sum_{k=1}^{\infty }\gamma _{j}^{\omega }\gamma _{k}^{\omega
}e_{k}(\Gamma e_{j},e_{k})_{H}=Fe_{j},\mbox{ }j\geq 1, 
\]%
which multiplied scalarly by $e_{j}$ yields 
\[
2\left\langle A_{\omega }e_{j},P_{\omega }e_{j}\right\rangle _{V^{\prime
},V}+(\gamma _{j}^{\omega })^{2}(\Gamma e_{j},e_{j})_{H}=(Fe_{j},e_{j})_{H},%
\mbox{ }j\geq 1. 
\]%
Finally, taking into account that $(Fe_{j},e_{j})_{H}=(C_{1}^{\ast
}C_{1}e_{j},e_{j})_{H}=\left\Vert C_{1}e_{j}\right\Vert _{H}^{2}>0,$ we
obtain the equation 
\[
(\gamma _{j}^{\omega })^{2}g_{jj}+2\gamma _{j}^{\omega }a_{jj}^{\omega
}-\left\Vert C_{1}e_{j}\right\Vert _{H}^{2}=0,\mbox{ }j\geq 1, 
\]%
where $g_{jj}:=(\Gamma e_{j},e_{j})_{H}$ and $a_{jj}^{\omega }:=\left\langle
A_{\omega }e_{j},e_{j}\right\rangle _{V^{\prime },V}.$ This equation has the
positive solution 
\begin{equation}
\gamma _{j}^{\omega }=\frac{-a_{jj}^{\omega }+\sqrt{(a_{jj}^{\omega
})^{2}+g_{jj}\left\Vert C_{1}e_{j}\right\Vert _{H}^{2}}}{g_{jj}}.
\label{55-7}
\end{equation}%
Since $g_{jj}\geq g_{0}>0$ by hypothesis (\ref{52}) and $a_{jj}^{\omega
}\geq \omega $, it follows that 
\begin{equation}
0<\gamma _{j}^{\omega }\leq \frac{\left\Vert C_{1}e_{j}\right\Vert _{H}}{%
\sqrt{g_{jj}}}\leq \frac{\left\Vert C_{1}e_{j}\right\Vert _{H}}{\sqrt{g_{0}}}%
\leq \frac{1}{\sqrt{g_{0}}}\left\Vert C_{1}\right\Vert _{\mathcal{H}}.
\label{55-6}
\end{equation}%
Then,%
\[
\left\Vert P_{\omega }\right\Vert _{\mathcal{H}}^{2}=\sum\limits_{j=1}^{%
\infty }\left\Vert P_{\omega }e_{j}\right\Vert
_{H}^{2}=\sum\limits_{j=1}^{\infty }|\gamma _{j}^{\omega }|^{2}\left\Vert
e_{j}\right\Vert _{H}^{2}\leq \sum\limits_{j=1}^{\infty }|\gamma
_{j}^{\omega }|^{2}\leq \frac{1}{g_{0}}\left\Vert C_{1}\right\Vert _{%
\mathcal{H}}^{2}. 
\]%
Moreover, recalling that $C_{1}\in (V^{\prime },H)_{H.S.}$ and $\mathcal{V}%
:=(V^{\prime },H)_{H.S.}\cap (H,V)_{H.S.}$ we have 
\begin{eqnarray}
&&\left\Vert P_{\omega }\right\Vert _{(V^{\prime
},H)_{H.S.}}^{2}=\sum\limits_{j=1}^{\infty }\left\Vert P_{\omega }\widetilde{%
\widetilde{e_{j}}}\right\Vert _{H}^{2}=\sum\limits_{j=1}^{\infty }|\gamma
_{j}^{\omega }|^{2}\rho _{j}^{2}\left\Vert e_{j}\right\Vert _{H}^{2}\leq 
\frac{1}{g_{0}}\sum\limits_{j=1}^{\infty }(C_{1}e_{j},C_{1}e_{j})_{H}\rho
_{j}^{2}  \label{55-8} \\
&=&\frac{1}{g_{0}}\sum\limits_{j=1}^{\infty }\left( C_{1}\widetilde{%
\widetilde{e_{j}}},C_{1}\widetilde{\widetilde{e_{j}}}\right) _{H}=\frac{1}{%
g_{0}}\left\Vert C_{1}\right\Vert _{(V^{\prime },H)_{H.S.}}^{2}  \nonumber
\end{eqnarray}%
and 
\begin{equation}
\left\Vert P_{\omega }\right\Vert
_{(H,V)_{H.S.}}^{2}=\sum\limits_{j=1}^{\infty }\left\Vert Pe_{j}\right\Vert
_{V}^{2}=\sum\limits_{j=1}^{\infty }|\gamma _{j}^{\omega }|^{2}\left\Vert
e_{j}\right\Vert _{V}^{2}\leq \frac{1}{g_{0}}\sum\limits_{j=1}^{\infty
}\left\Vert C_{1}e_{j}\right\Vert _{H}^{2}\rho _{j}^{2}\leq \frac{1}{g_{0}}%
\left\Vert C_{1}\right\Vert _{(V^{\prime },H)_{H.S.}}^{2},  \label{55-9}
\end{equation}%
thus 
\begin{equation}
\left\Vert P_{\omega }\right\Vert _{\mathcal{V}}\leq \frac{2}{\sqrt{g_{0}}}%
\left\Vert C_{1}\right\Vert _{(V^{\prime },H)_{H.S.}}.  \label{55-3}
\end{equation}%
Then, by (\ref{24}) we deduce that 
\[
\left\Vert \mathcal{B}(P_{\omega })\right\Vert _{\mathcal{H}}\leq \left\Vert
\Gamma \right\Vert _{\mathcal{H}}\left\Vert P_{\omega }\right\Vert _{%
\mathcal{H}}^{2}\leq C 
\]%
and so, by the equation (\ref{54}), $\mathcal{L}_{\omega }(P_{\omega })=F-%
\mathcal{B}(P_{\omega })$ and it follows that 
\[
\left\Vert \mathcal{L}_{\omega }(P_{\omega })\right\Vert _{\mathcal{H}}\leq
C, 
\]%
where $C$ denotes several positive constants independent of $\omega .$

By (\ref{55-3}) it is obvious that the operator $\mathcal{L}$ satisfies 
\begin{equation}
\left\Vert \mathcal{L}(P_{\omega })\right\Vert _{\mathcal{V}^{\prime }}\leq
C.  \label{55-4}
\end{equation}%
By the above boundedness we deduce that there exist $P,$ $L,$ $L_{1}$ and $%
\eta $ in $\mathcal{H}$ such that 
\begin{equation}
P_{\omega }\rightarrow P,\mbox{ weakly in }\mathcal{V},  \label{55-1}
\end{equation}%
\[
\mathcal{B}(P_{\omega })\rightarrow \eta ,\mbox{ weakly in }\mathcal{H}, 
\]%
\[
\mathcal{L}_{\omega }(P_{\omega })\rightarrow L:=F-\eta \mbox{, weakly in }%
\mathcal{H}\mbox{ and }\mathcal{V}^{\prime }, 
\]%
\[
\mathcal{L}(P_{\omega })\rightarrow L_{1},\mbox{ weakly in }\mathcal{V}%
^{\prime }. 
\]

Moreover, by (\ref{55-6}), there exist $\gamma _{j}\geq 0$ such that $\gamma
_{j}^{\omega }\rightarrow \gamma _{j}$ as $\omega \rightarrow 0,$ for all $%
j, $ and it follows that $Pe_{j}=\gamma _{j}e_{j}.$ Indeed, for $Z\in 
\mathcal{H},$ we have 
\[
\left( P_{\omega },Z\right) _{\mathcal{H}}=\sum\limits_{j=1}^{\infty }\left(
P_{\omega }e_{j},Ze_{j}\right) _{H}=\sum\limits_{j=1}^{\infty }\left( \gamma
_{j}^{\omega }e_{j},Ze_{j}\right) _{H}\rightarrow \sum\limits_{j=1}^{\infty
}\left( \gamma _{j}e_{j},Ze_{j}\right) _{H}=\left( P,Z\right) _{\mathcal{H}%
}. 
\]%
Then we have 
\begin{eqnarray}
(\mathcal{B}(P_{\omega }),P_{\omega })_{\mathcal{H}}
&=&\sum\limits_{j=1}^{\infty }(P_{\omega }\Gamma P_{\omega }e_{j},P_{\omega
}e_{j})_{H}=\sum\limits_{j=1}^{\infty }(\gamma _{j}^{\omega })^{2}g_{jj}
\label{55-5} \\
&\rightarrow &\sum\limits_{j=1}^{\infty }\gamma
_{j}^{2}g_{jj}=\sum\limits_{j=1}^{\infty }(P\Gamma Pe_{j},Pe_{j})_{H}=(%
\mathcal{B}(P),P)_{\mathcal{H}},  \nonumber
\end{eqnarray}%
hence $\eta =\mathcal{B}(P).$

We have to show that $L=\mathcal{L}(P).$ To this end we recall that if $%
\mathcal{L}$ is maximal monotone, $P_{\omega }\rightarrow P$ weakly in $%
\mathcal{V},$ $\mathcal{L(}P_{\omega })\rightarrow L_{1}$ weakly in $%
\mathcal{V}^{\prime }$ and if 
\begin{equation}
\limsup\limits_{\omega \rightarrow 0}\left\langle \mathcal{L}(P_{\omega
}),P_{\omega }\right\rangle _{\mathcal{V}^{\prime },\mathcal{V}}\leq
\left\langle L_{1},P\right\rangle _{\mathcal{V}^{\prime },\mathcal{V}}
\label{55-0}
\end{equation}%
then $L_{1}=\mathcal{L}(P)$ (see \cite{VB-2010}, p.41, Corollary 2.4).

It remains to prove that $\mathcal{L}$ is maximal monotone and to show (\ref%
{55-0}). First, we have $\left\langle \mathcal{L}P,P\right\rangle
_{V^{\prime },V}\geq 0.$ Then, we have to show that equation 
\begin{equation}
P+\lambda \mathcal{L}(P)=G\in \mathcal{H}  \label{55-2}
\end{equation}%
has a solution $P\in D(\mathcal{L}).$ To prove this we start from the
equation 
\[
P_{\omega }+\lambda \mathcal{L}_{\omega }(P_{\omega })=G, 
\]%
or equivalently%
\[
\sum_{j=1}^{\infty }y_{j}\gamma _{j}^{\omega }e_{j}+\lambda
\sum_{j=1}^{\infty }y_{j}(A_{\omega }^{\ast }P_{\omega }+P_{\omega
}A_{\omega })e_{j}=\sum_{j=1}^{\infty }y_{j}Ge_{j},\mbox{ for all }y\in H. 
\]%
After some algebra we obtain 
\[
\gamma _{j}^{\omega }=\left( 1+2\lambda \left\langle A_{\omega
}e_{j},e_{j}\right\rangle _{V^{\prime },V}\right) ^{-1}(Ge_{j},e_{j})_{H}, 
\]%
which leads, at limit as $\omega \rightarrow 0,$ to 
\[
\gamma _{j}=\left( 1+2\lambda \left\langle Ae_{j},e_{j}\right\rangle
_{V^{\prime },V}\right) ^{-1}(Ge_{j},e_{j})_{H}. 
\]%
This proves, by a similar calculation, that $P$ defined by $Pe_{j}=\gamma
_{j}e_{j}$ solves equation (\ref{55-2}), hence $\mathcal{L}$ is maximal
monotone on $\mathcal{H}.$

By (\ref{54}) and (\ref{55-5}) we have 
\begin{equation}
\left\langle \mathcal{L}_{\omega }(P_{\omega }),P_{\omega }\right\rangle _{%
\mathcal{V}^{\prime },\mathcal{V}}=(F-\mathcal{B}(P_{\omega }),P_{\omega })_{%
\mathcal{H}}\rightarrow (F,P)_{\mathcal{H}}-(\mathcal{B}(P),P)_{\mathcal{H}}
\label{55-0-2}
\end{equation}%
and write the left-hand side term as 
\[
\left\langle \mathcal{L}_{\omega }(P_{\omega }),P_{\omega }\right\rangle _{%
\mathcal{V}^{\prime },\mathcal{V}}=\left\langle \mathcal{L}_{\omega
}(P_{\omega })-\mathcal{L}(P_{\omega }),P_{\omega }\right\rangle _{\mathcal{V%
}^{\prime },\mathcal{V}}+\left\langle \mathcal{L}(P_{\omega }),P_{\omega
}\right\rangle _{\mathcal{V}^{\prime },\mathcal{V}}. 
\]%
Then, taking into account that $P_{\omega }\in \mathcal{V}=(V^{\prime
},H)_{H.S.}\cap (H,V)_{H.S.}$ we have 
\begin{eqnarray*}
&&\left\langle \mathcal{L}_{\omega }(P_{\omega })-\mathcal{L}(P_{\omega
}),P_{\omega }\right\rangle _{\mathcal{V}^{\prime },\mathcal{V}}=\omega
\sum\limits_{j=1}^{\infty }\left\langle JP_{\omega }+P_{\omega }J,P_{\omega
}\right\rangle _{\mathcal{V}^{\prime },\mathcal{V}} \\
&=&\omega \sum\limits_{j=1}^{\infty }\left\langle JP_{\omega },P_{\omega
}\right\rangle _{(H,V^{\prime })_{H.S.},(H,V)_{H.S.}}+\omega
\sum\limits_{j=1}^{\infty }\left\langle P_{\omega }J,P_{\omega
}\right\rangle _{(V,H)_{H.S.},(V^{\prime },H)_{H.S.}} \\
&=&\omega \sum\limits_{j=1}^{\infty }\left\langle JP_{\omega
}e_{j},P_{\omega }e_{j}\right\rangle _{V^{\prime },V}+\omega
\sum\limits_{j=1}^{\infty }\left( P_{\omega }J\widetilde{e_{j}},P_{\omega }%
\widetilde{\widetilde{e_{j}}}\right) _{H} \\
&=&\omega \sum\limits_{j=1}^{\infty }\left\Vert P_{\omega }e_{j}\right\Vert
_{V}^{2}+\omega \sum\limits_{j=1}^{\infty }\left\langle J\widetilde{e_{j}}%
,\rho _{j}P_{\omega }^{2}e_{j}\right\rangle _{V^{\prime },V} \\
&\leq &\omega \left\Vert P_{\omega }\right\Vert _{(H,V)_{H.S.}}^{2}+\omega
\sum\limits_{j=1}^{\infty }\left\langle J\widetilde{e_{j}},\rho
_{j}^{2}(\gamma _{j}^{\omega })^{2}\widetilde{e_{j}}\right\rangle
_{V^{\prime },V}\leq \frac{\omega }{g_{0}}\left\Vert C_{1}\right\Vert
_{(V^{\prime },H)_{H.S.}}^{2}+\omega \sum\limits_{j=1}^{\infty }\left(
\gamma _{j}^{\omega }\right) ^{2}\left\Vert \widetilde{e_{j}}\right\Vert
_{V}^{2}\rho _{j}^{2} \\
&\leq &\frac{\omega }{g_{0}}\left\Vert C_{1}\right\Vert _{(V^{\prime
},H)_{H.S.}}^{2}+\frac{\omega }{g_{0}}\sum\limits_{j=1}^{\infty
}(C_{1}e_{j},C_{1}e_{j})_{H}\rho _{j}^{2}=\frac{2\omega }{g_{0}}\left\Vert
C_{1}\right\Vert _{(V^{\prime },H)_{H.S.}}^{2}.
\end{eqnarray*}%
In the last inequalities we used (\ref{55-9}). Therefore,%
\[
\lim\limits_{\omega \rightarrow 0}\left\langle \mathcal{L}_{\omega
}(P_{\omega })-\mathcal{L}(P_{\omega }),P_{\omega }\right\rangle _{\mathcal{V%
}^{\prime },\mathcal{V}}=0. 
\]%
Writing now 
\[
\left\langle \mathcal{L}_{\omega }(P_{\omega }),Z\right\rangle _{\mathcal{V}%
^{\prime },\mathcal{V}}=\left\langle \mathcal{L}_{\omega }(P_{\omega })-%
\mathcal{L}(P_{\omega }),Z\right\rangle _{\mathcal{V}^{\prime },\mathcal{V}%
}+\left\langle \mathcal{L}(P_{\omega }),Z\right\rangle _{\mathcal{V}^{\prime
},\mathcal{V}},\mbox{ for all }Z\in \mathcal{V}, 
\]%
we deduce that 
\[
\lim_{\omega \rightarrow 0}\left\langle \mathcal{L}(P_{\omega
}),Z\right\rangle _{\mathcal{V}^{\prime },\mathcal{V}}=\lim_{\omega
\rightarrow 0}\left\langle \mathcal{L}_{\omega }(P_{\omega }),Z\right\rangle
_{\mathcal{V}^{\prime },\mathcal{V}} 
\]%
hence $\mathcal{L}(P_{\omega })\rightarrow L_{1}=F-\mathcal{B}(P)$ weakly in 
$\mathcal{V}^{\prime }.$

We have 
\begin{eqnarray}
&&\limsup\limits_{\omega \rightarrow 0}\left\langle \mathcal{L}(P_{\omega
}),P_{\omega }\right\rangle _{\mathcal{V}^{\prime },\mathcal{V}%
}=\limsup\limits_{\omega \rightarrow 0}\left\{ \left\langle \mathcal{L}%
_{\omega }(P_{\omega }),P_{\omega }\right\rangle _{\mathcal{V}^{\prime },%
\mathcal{V}}-\left\langle \mathcal{L}_{\omega }(P_{\omega })-\mathcal{L}%
(P_{\omega }),P_{\omega }\right\rangle _{\mathcal{V}^{\prime },\mathcal{V}%
}\right\}  \nonumber \\
&=&\limsup\limits_{\omega \rightarrow 0}(F-\mathcal{B}(P_{\omega
}),P_{\omega })_{\mathcal{H}}=\lim_{\omega \rightarrow 0}(F,P_{\omega })_{%
\mathcal{H}}-\lim\limits_{\omega \rightarrow 0}(\mathcal{B}(P_{\omega
}),P_{\omega })_{\mathcal{H}}  \label{58}
\end{eqnarray}%
and using (\ref{55-0-2}) we deduce%
\[
\limsup\limits_{\omega \rightarrow 0}(\mathcal{L}(P_{\omega }),P_{\omega })_{%
\mathcal{H}}\leq (F,P)_{\mathcal{H}}-(\mathcal{B}(P),P)_{\mathcal{H}}=(F-%
\mathcal{B}(P),P)_{\mathcal{H}}=(L,P)_{\mathcal{H}}. 
\]%
We conclude that $\mathcal{L}(P)=L_{1}=F-\mathcal{B}(P)=L,$ according to the
previous mentioned result in \cite{VB-2010}.

Now, we have all ingredients to pass to the limit in (\ref{54}) as $\omega
\rightarrow 0$ and obtain that $\mathcal{L}(P)+\mathcal{B}(P)=F,$ that is $P$
is the solution to (\ref{53}), as claimed.

Finally, we prove that the solution to (\ref{53}) is unique.

First of all we observe that the eigenvalues $\gamma _{j}$ of a solution $P$
to (\ref{53}) can be obtained by the same calculus as that leading to (\ref%
{55-7}), where $a_{ij}^{\omega }$ is replaced by $a_{ij}=\left\langle
Ae_{j},e_{j}\right\rangle _{V^{\prime },V}$ and that these eigenvalues $%
\gamma _{j}$ are positive since $g_{jj}\geq g_{0}>0$ and $\left\Vert
C_{1}e_{j}\right\Vert _{H}>0,$ under the assumption $\ker C_{1}=\{0\}.$

Let us consider two solutions $P_{1}$ and $P_{2}$ to (\ref{53}), that is 
\begin{equation}
A^{\ast }P_{i}y+P_{i}Ay+P_{i}\Gamma P_{i}y=C_{1}^{\ast }C_{1}y,%
\mbox{ for
all }y\in H,\mbox{ }i=1,2,  \label{58-1}
\end{equation}%
and recall that $P_{i}$ have the representations%
\[
P_{1}y=\sum_{j=1}^{\infty }\gamma _{j}^{1}y_{j}e_{j},\mbox{ }%
P_{2}y=\sum_{j=1}^{\infty }\gamma _{j}^{2}y_{j}e_{j} 
\]%
where $\{e_{j}\}_{j}$ is the basis corresponding to the canonical
isomorphism $J.$ Moreover, the eigenvalues $\gamma _{j}^{1}$ and $\gamma
_{j}^{2}$ corresponding to $P_{1}$ and $P_{2},$ respectively, are positive.
Denote 
\begin{equation}
\gamma _{j}:=\gamma _{j}^{1}-\gamma _{j}^{2},\mbox{ }P:=P_{1}-P_{2}=%
\sum_{j=1}^{\infty }\gamma _{j}y_{j}e_{j}  \label{58-2}
\end{equation}%
and subtract equations (\ref{58-1}) corresponding to $P_{1}$ and $P_{2}$,
getting%
\[
A^{\ast }Py+PAy+P\Gamma P_{1}y+P_{2}\Gamma Py=0,\mbox{ for all }y\in H. 
\]%
Replacing $y$ by its expansion with respect to the same basis, we have 
\[
\sum_{j=1}^{\infty }y_{j}(A^{\ast }P+PA)e_{j}+\sum_{j=1}^{\infty
}y_{j}\sum_{i=1}^{\infty }(\gamma _{j}\gamma _{i}^{2}+\gamma _{i}\gamma
_{j}^{1})g_{ji}e_{i}=0,\mbox{ for all }y. 
\]%
We multiply scalarly by $e_{j}$ and finally obtain 
\[
2\sum_{j=1}^{\infty }y_{j}\left\langle Ae_{j},Pe_{j}\right\rangle
_{V^{\prime },V}+\sum_{j=1}^{\infty }y_{j}\gamma _{j}(\gamma _{j}^{1}+\gamma
_{j}^{2})g_{jj}=0 
\]%
which holds for all $y,$ hence we deduce using $P$ written in (\ref{58-2})
that%
\[
\gamma _{j}\left( 2a_{jj}+(\gamma _{j}^{1}+\gamma _{j}^{2})g_{jj}\right) =0,%
\mbox{ }j\geq 1. 
\]%
Since $2a_{jj}+(\gamma _{j}^{1}+\gamma _{j}^{2})g_{jj}>0$, this implies that 
$\gamma _{j}=0,$ which shows that $P_{1}=P_{2},$ ending thus the
proof.\hfill $\square $

\section{Solving the $H^{\infty }$-control problem}

\setcounter{equation}{0}

In this section we give the result stating that the $H^{\infty }$-control
problem associated to system (\ref{1})-(\ref{2}) has a solution.

\begin{theorem}
\label{T4.1}Let $\gamma >0$. Let the assumptions of Theorem 2.1 and
Corollary 2.2 , or Theorem 3.1, respectively hold and let $w\in L^{2}(%
\mathbb{R}_{+};W).$ Assume still that $-A-B_{2}B_{2}^{\ast }P$ generates an
exponentially stable semigroup on $H,$ where $P$ is the solution to the
Riccati equation (\ref{5}). Then, the feedback operator $F=-B_{2}^{\ast }P$
solves the $H^{\infty }$-control problem associated to system (\ref{1})-(\ref%
{2}).
\end{theorem}

\noindent \textbf{Proof}. Under the hypotheses of the theorems indicated
above it follows that the Riccati equation (\ref{5}) has a unique solution.
Let us set $F=-B_{2}^{\ast }P$ and introduce it in equation (\ref{1}),
obtaining 
\begin{equation}
y^{\prime }(t)=-Ay(t)-B_{2}B_{2}^{\ast }Py(t)+B_{1}w(t),\mbox{ }t>0,
\label{70}
\end{equation}%
with $y_{0}=0.$ We have to show that the solution to this equation belongs
to $L^{2}(\mathbb{R}_{+};H)$ and that $z(t)=\Phi (y(t),u(t))$ verifies
inequality (\ref{3}). Denote $\widetilde{A}=-A-B_{2}B_{2}^{\ast }P,$ write
the solution as 
\[
y(t)=\int_{0}^{t}e^{\widetilde{A}(t-s)}B_{1}w(s)ds 
\]%
and calculate 
\[
\left\Vert y(t)\right\Vert _{H}\leq \int_{0}^{t}\left\Vert e^{\widetilde{A}%
(t-s)}B_{1}w(s)\right\Vert _{H}ds\leq \int_{0}^{t}\left\Vert e^{-\alpha
(t-s)}B_{1}w(s)\right\Vert _{H}ds\leq C\int_{0}^{t}e^{-\alpha
(t-s)}\left\Vert w(s)\right\Vert _{H}ds, 
\]%
since $\widetilde{A}$ generates an exponentially stable semigroup on $H$
(see (\ref{2-0})). Then, we integrate and apply the Young's inequality for
convolution,%
\begin{eqnarray*}
\int_{0}^{\infty }\left\Vert y(t)\right\Vert _{H}^{2}dt &\leq
&C\int_{0}^{\infty }\left( \int_{0}^{t}e^{-\alpha (t-s)}\left\Vert
w(s)\right\Vert _{H}ds\right) ^{2}dt\leq C\left( \int_{0}^{\infty
}e^{-\alpha t}dt\right) ^{2}\int_{0}^{\infty }\left\Vert w(t)\right\Vert
^{2}dt \\
&\leq &C\left\Vert w\right\Vert _{L^{2}(\mathbb{R}_{+}:H)}^{2}<\infty .
\end{eqnarray*}%
We pass now to the verification of (\ref{3}). We multiply (\ref{70}) by $%
Py(t),$ getting 
\[
\frac{1}{2}\frac{d}{dt}(Py(t),y(t))_{H}+\left\langle
Ay(t),Py(t)\right\rangle _{V^{\prime },V}=(B_{1}w(t),Py(t))_{H}-\left\Vert
B_{2}^{\ast }Py(t)\right\Vert _{H}^{2} 
\]%
which integrated with respect to $t,$ yields%
\begin{equation}
\int_{0}^{\infty }\left( 2\left\langle Ay(t),Py(t)\right\rangle _{V^{\prime
},V}-2(B_{1}w(t),Py(t))_{H}+2\left\Vert B_{2}^{\ast }Py(t)\right\Vert
_{H}^{2}\right) dt=0.  \label{71}
\end{equation}%
The Riccati equation (\ref{5}) yields 
\[
2\left\langle Ay(t),Py(t)\right\rangle _{V^{\prime },V}+\left\Vert
B_{2}^{\ast }Py(t)\right\Vert _{H}^{2}-\gamma ^{-2}\left\Vert B_{1}^{\ast
}Py(t)\right\Vert _{H}^{2}=\left\Vert C_{1}y(t)\right\Vert _{H}^{2} 
\]%
which replaced in (\ref{71}) gives 
\[
\int_{0}^{\infty }(\left\Vert C_{1}y(t)\right\Vert _{H}^{2}+\left\Vert
B_{2}^{\ast }Py(t)\right\Vert _{H}^{2}+\gamma ^{-2}\left\Vert B_{1}^{\ast
}Py(t)\right\Vert _{H}^{2}-2(B_{1}w(t),Py(t))_{H})dt=0. 
\]%
This can be still written%
\begin{eqnarray*}
&&\int_{0}^{\infty }(\left\Vert C_{1}y(t)\right\Vert _{H}^{2}+\left\Vert
B_{2}^{\ast }Py(t)\right\Vert _{H}^{2})dt+\int_{0}^{\infty }\gamma
^{2}\left\Vert w(t)-\gamma ^{-2}B_{1}^{\ast }Py(t)\right\Vert _{W}^{2}dt \\
&=&\int_{0}^{\infty }\gamma ^{2}\left\Vert w(t)\right\Vert _{W}^{2}dt.
\end{eqnarray*}%
By denoting $\widetilde{w}(t)=w(t)-\gamma ^{-2}B_{1}^{\ast }Py(t)$ we get 
\[
\int_{0}^{\infty }(\left\Vert C_{1}y(t)\right\Vert _{H}^{2}+\left\Vert
B_{2}^{\ast }Py(t)\right\Vert _{H}^{2})dt\leq \gamma ^{2}\int_{0}^{\infty
}(\left\Vert w(t)\right\Vert _{W}^{2}-\left\Vert \widetilde{w}(t)\right\Vert
_{W}^{2}dt. 
\]%
Now, we assert that there exists $0<\alpha <1$ such that $\left\Vert 
\widetilde{w}\right\Vert _{L^{2}(\mathbb{R}_{+}:W)}^{2}\geq \alpha
\left\Vert w\right\Vert _{L^{2}(\mathbb{R}_{+};W)}^{2},$ for all $w\in L^{2}(%
\mathbb{R}_{+};H).$ This result is proved in Lemma 3.5 in \cite{GM-23-Hinf}.
Then, it follows that%
\begin{eqnarray*}
\left\Vert z\right\Vert _{L^{2}(\mathbb{R}_{+};Z)}^{2} &=&\left\Vert \Phi
(y,u)\right\Vert _{L^{2}(\mathbb{R}_{+};Z)}^{2}\leq \int_{0}^{\infty
}(\left\Vert C_{1}y(t)\right\Vert _{H}^{2}+\left\Vert B_{2}^{\ast
}Py(t)\right\Vert _{H}^{2})dt \\
&\leq &\gamma ^{2}(1-\alpha )\left\Vert w\right\Vert _{L^{2}(\mathbb{R}%
_{+};W)}^{2}<\gamma ^{2}\left\Vert w\right\Vert _{L^{2}(\mathbb{R}%
_{+};W)}^{2},
\end{eqnarray*}%
by (\ref{2-0}), as claimed. This ends the proof.\hfill $\square $

\section{An example}

\setcounter{equation}{0}

We shall give an example for applying Theorem \ref{T2.1} to a parabolic
equation with a singular term of Hardy type.

Let $\Omega $ be an open bounded subset of $\mathbb{R}^{N},$ with the
boundary $\Gamma =\partial \Omega $ sufficiently smooth. We consider the
following singular system%
\begin{equation}
y_{t}-\Delta y-\frac{\lambda y}{|x|^{2}}=B_{1}w+B_{2}u,\mbox{ in }(0,\infty
)\times \Omega ,  \label{21-0}
\end{equation}%
\begin{equation}
y=0,\mbox{ on }(0,\infty )\times \partial \Omega ,  \label{21-1}
\end{equation}%
\begin{equation}
y(0)=0,\mbox{ in }\Omega ,  \label{21-2}
\end{equation}%
\begin{equation}
z=C_{1}y,\mbox{ in }(0,\infty )\times \Omega ,  \label{21-3}
\end{equation}%
where $\left\vert \cdot \right\vert $ denotes the Euclidian norm in $\mathbb{%
R}^{N}$ and $\lambda $ and $N$ will be specified a few lines below.

Let $k\in L^{2}(\Omega \times \Omega ).$ We define the operators 
\begin{eqnarray}
B_{1} &:&L^{2}(\Omega )\rightarrow H_{0}^{1}(\Omega ),\mbox{ }%
B_{1}w=(-\Delta )^{-1}w,\mbox{ \ }  \label{21-7} \\
B_{2} &:&L^{2}(\Omega )\rightarrow L^{2}(\Omega ),\mbox{ }%
(B_{2}u)(x)=\int_{\Omega }k(x,\xi )u(\xi )d\xi ,  \nonumber \\
C_{1} &:&L^{2}(\Omega )\rightarrow H_{0}^{1}(\Omega ),\mbox{ }%
C_{1}y=(-\Delta )^{-1}y.  \nonumber
\end{eqnarray}%
In this case, we consider the variational triplet $H_{0}^{1}(\Omega )\subset
L^{2}(\Omega )\subset H^{-1}(\Omega ),$ define 
\[
-\Delta :H_{0}^{1}(\Omega )\rightarrow H^{-1}(\Omega ),\mbox{ }\left\langle
-\Delta y,\psi \right\rangle _{V^{\prime },V}=\int_{\Omega }\nabla y\cdot
\nabla \psi dx 
\]%
and choose%
\begin{equation}
H=W=Z=L^{2}(\Omega ),\mbox{ }V=H_{0}^{1}(\Omega ),\mbox{ }U=L^{2}(\Omega ).
\label{21-6}
\end{equation}%
We introduce the self-adjoint operator 
\[
A:D(A)\subset L^{2}(\Omega )\rightarrow L^{2}(\Omega ),\mbox{ }Ay=-\Delta y-%
\frac{\lambda y}{|x|^{2}}, 
\]%
with $D(A)=\{y\in H_{0}^{1}(\Omega );$ $Ay\in L^{2}(\Omega )\}.$ Then,
equation (\ref{21-0}) can be equivalently written%
\begin{equation}
y^{\prime }(t)+Ay(t)=B_{1}w(t)+B_{2}u(t),\mbox{ }t\geq 0.  \label{24-1}
\end{equation}%
In order to apply Theorem \ref{T2.1} and to prove the existence of the
Riccati equation (\ref{6}) in this case, we have to check conditions (\ref{8}%
), (\ref{12}) and (\ref{13}) for the operator $A$ and for the operators $%
B_{1},$ $B_{2},$ $C_{1}.$

Now, we discuss the issues possibly raised by the term $\frac{\lambda y}{%
|x|^{2}}$ in some situations$.$

If $0\in \Omega ,$ this term is a singular contribution of the Hardy type
which requires the application of the Hardy inequality 
\begin{equation}
\int_{\Omega }\left\vert \nabla y(x)\right\vert ^{2}dx\geq H_{N}\int_{\Omega
}\frac{\left\vert y(x)\right\vert ^{2}}{\left\vert x\right\vert ^{2}}dx,%
\mbox{ for all }y\in H_{0}^{1}(\Omega ),  \label{HN}
\end{equation}%
where $H_{N}=\frac{(N-2)^{2}}{4}$ is optimal Hardy constant (see \cite%
{Brezis-Vazquez}, p. 452, Theorem 4.1) and (\ref{HN}) takes place for $N\geq
3.$ We consider the case $0<\lambda <H_{N}.$

First, we verify (\ref{8}). In this case the operator $A$ is $m$-accretive
on $L^{2}(\Omega )$ (see \cite{GM-23-Hinf}, Section 4) and coercive$.$ To
this end we shall use several times the Hardy inequality (\ref{HN}) which
also ensures that $\frac{y}{x}\in L^{2}(\Omega )$ if $y\in H_{0}^{1}(\Omega
).$ We have 
\begin{eqnarray*}
\left\langle Ay,y\right\rangle _{V^{\prime },V} &=&\int_{\Omega }\left\vert
\nabla y\right\vert ^{2}dx-\lambda \int_{\Omega }\frac{\left\vert
y\right\vert ^{2}}{\left\vert x\right\vert ^{2}}dx\geq \left( 1-\frac{%
\lambda }{H_{N}}\right) \int_{\Omega }\left\vert \nabla y\right\vert ^{2}dx
\\
&\geq &\frac{1}{2}\left( 1-\frac{\lambda }{H_{N}}\right) \left\Vert \nabla
y\right\Vert _{H}^{2}+\frac{H_{N}}{2}\left( 1-\frac{\lambda }{H_{N}}\right)
\left\Vert \frac{y}{x}\right\Vert _{H}^{2}.
\end{eqnarray*}%
If $0\notin \Omega $ the problem is studied for $N\geq 1$ and for $\lambda
\geq 0$ small enough (if different of $0),$ case in which standard
calculations are done, so that we shall not treat it in detail.

It remains to check hypotheses (\ref{12}) and (\ref{13}).

We see that 
\[
B_{1}\in L(H,H),\mbox{ }B_{2}\in L(H,H),\mbox{ }C_{1}\in L(H,H)
\]%
and $B_{1}^{\ast }$ and $B_{2}^{\ast }$ are defined by 
\begin{eqnarray}
B_{1}^{\ast } &:&H\rightarrow V,\mbox{ }B_{1}^{\ast }v=(-\Delta )^{-1}v,%
\mbox{ for }v\in H,  \label{22-3} \\
B_{2}^{\ast } &:&H\rightarrow H,\mbox{ }(B_{2}^{\ast }v)(x)=\int_{\Omega
}k(\xi ,x)v(\xi )d\xi ,\mbox{ for }v\in H,  \nonumber \\
C_{1}^{\ast } &:&H\rightarrow V,\mbox{ }C_{1}^{\ast }v=(-\Delta )^{-1}v,%
\mbox{ for }v\in H.  \nonumber
\end{eqnarray}%
Hence 
\begin{eqnarray*}
B_{1}B_{1}^{\ast } &:&H\rightarrow H,\mbox{ }B_{1}B_{1}^{\ast }v=(-\Delta
)^{-2}v, \\
B_{2}B_{2}^{\ast } &:&H\rightarrow H,\mbox{ }(B_{2}B_{2}^{\ast
}v)(x)=\int_{\Omega }\int_{\Omega }k(x,\xi )k(\xi ^{\prime },\xi )v(\xi
^{\prime })d\xi ^{\prime }d\xi , \\
C_{1}^{\ast }C_{1} &:&H\rightarrow H,\mbox{ }C_{1}C_{1}^{\ast }v=(-\Delta
)^{-2}v.
\end{eqnarray*}%
We show that $\Gamma =B_{2}B_{2}^{\ast }-\gamma ^{-2}B_{1}B_{1}^{\ast }$ and 
$C_{1}^{\ast }C_{1}$ belong to $\mathcal{H}=(H,H)_{H.S.}$. We choose the
basis $\{e_{j}\}_{j}$ generated by the canonical isomorphism $-\Delta
:H_{0}^{1}(\Omega )\rightarrow H^{-1}(\Omega ),$ 
\[
-\Delta e_{j}=\lambda _{j}e_{j}\mbox{ in }\Omega ,\mbox{ }e_{j}=0\mbox{ on }%
\partial \Omega .
\]%
We have 
\[
\left\Vert B_{1}B_{1}^{\ast }\right\Vert _{\mathcal{H}}=\sum\limits_{j=1}^{%
\infty }\left\Vert (-\Delta )^{-1}(-\Delta )^{-1}e_{j}\right\Vert
_{H}^{2}=\sum\limits_{j=1}^{\infty }\lambda _{j}^{-2},
\]%
and recall that, by the Weyl's law characterizing the asymptotic
distribution of the eigenvalues for the Laplacian on a bounded domain, there
is a positive constant $c_{1}^{\Omega }$ depending on $\Omega ,$ such that 
\begin{equation}
c_{1}^{\Omega }j^{2/N}<\lambda _{j},\mbox{ for all }j\geq 1,  \label{59}
\end{equation}%
(see e.g., \cite{Li-Yau-1983}). Hence%
\[
\left\Vert B_{1}B_{1}^{\ast }\right\Vert _{\mathcal{H}}\leq \frac{1}{%
(c_{1}^{\Omega })^{2}}\sum\limits_{j=1}^{\infty }j^{-4/N}<\infty ,%
\mbox{ for 
}1\leq N<4.
\]%
Similarly, we obtain $\left\Vert C_{1}C_{1}^{\ast }\right\Vert _{\mathcal{H}%
}<\infty .$ Thus, if $0<\lambda <H_{N},$ the above relations take place for $%
N=3,$ due to the dimension restriction impose by the presence of the Hardy
term. In the absence of the Hardy term, when $\lambda =0,$ the relations are
satisfied if $N\leq 3.$

Next, by the Parseval equality, 
\[
\left\Vert B_{2}\right\Vert _{\mathcal{H}}^{2}=\sum\limits_{j=1}^{\infty
}\int_{\Omega }\left( \int_{\Omega }k(x,\xi )e_{j}(\xi )d\xi \right)
^{2}dx=\int_{\Omega }\int_{\Omega }k^{2}(x,\xi )d\xi dx<\infty . 
\]%
Similarly $B_{2}^{\ast }\in \mathcal{H}$ and this implies that $%
B_{2}B_{2}^{\ast }\in \mathcal{H}$ by property (c) in Section 1.1. Thus, it
follows that $\Gamma \in \mathcal{H}$ and in the same way it is shown that $%
C_{1}^{\ast }C_{1}\in \mathcal{H}.$ Moreover, it is obvious that $%
C_{1}^{\ast }C_{1}\geq 0.$

Now, we shall indicate a sufficient condition involving the operators $B_{1}$
and $B_{2}$ which ensures that $\Gamma \geq 0.$

\begin{lemma}
\label{L4.1}Let $0<\lambda <H_{N}$, $N=3$ and assume the following
conditions:%
\begin{equation}
\sum\limits_{j=1}^{\infty }j^{-4/N}\leq (c_{1}^{\Omega })^{2}\gamma
^{2}\lambda _{1}^{B},  \label{64}
\end{equation}%
\begin{equation}
B_{2}B_{2}^{\ast }>0,\mbox{ }  \label{64-0}
\end{equation}%
where $\lambda _{1}^{B}$ is the first eigenvalue of the operator $%
B_{2}B_{2}^{\ast }$ and $c_{1}^{\Omega }$ is a positive constant depending
on $\Omega .$ Then, $\Gamma =B_{2}B_{2}^{\ast }-\gamma ^{-2}B_{1}B_{1}^{\ast
}\geq 0.$

If $\lambda =0$, the the result holds for all $N\leq 3.$
\end{lemma}

\noindent \textbf{Proof.} The condition $\Gamma \geq 0,$ that is $(\Gamma
y,y)_{H}\geq 0$ for all $y\in H,$ reduces to%
\[
(\Gamma y,y)=(B_{2}B_{2}^{\ast }y-\gamma ^{-2}B_{1}B_{1}^{\ast
}y,y)_{H}=\left\Vert B_{2}^{\ast }y\right\Vert _{H}^{2}-\gamma
^{-2}\left\Vert B_{1}^{\ast }y\right\Vert _{H}^{2}\geq 0,\mbox{ for all }%
y\in H,
\]%
and implies that%
\begin{equation}
\left\Vert (-\Delta )^{-1}y\right\Vert _{H}^{2}\leq \gamma ^{2}\left\Vert
B_{2}^{\ast }y\right\Vert _{H}^{2},\mbox{ for all }y\in H.  \label{60}
\end{equation}%
By (\ref{59}) we have 
\begin{equation}
\left\Vert (-\Delta )^{-1}y\right\Vert _{H}^{2}=\sum\limits_{j=1}^{\infty
}\lambda _{j}^{-2}(y,e_{j})_{H}^{2}\leq (c_{1}^{\Omega })^{-2}\left\Vert
y\right\Vert _{H}^{2}\sum\limits_{j=1}^{\infty }j^{-4/N}.  \label{61}
\end{equation}%
By comparison with (\ref{60}) it is sufficient to have 
\begin{equation}
(c_{1}^{\Omega })^{-2}\left\Vert y\right\Vert
_{H}^{2}\sum\limits_{j=1}^{\infty }j^{-4/N}\leq \gamma ^{2}\left\Vert
B_{2}^{\ast }y\right\Vert _{H}^{2}  \label{61-0}
\end{equation}%
or more, 
\begin{equation}
\sum\limits_{j=1}^{\infty }j^{-4/N}\leq \gamma ^{2}(c_{1}^{\Omega
})^{2}\inf_{y\in H,\mbox{ }y\neq 0}\frac{\left\Vert B_{2}^{\ast
}y\right\Vert _{H}^{2}}{\left\Vert y\right\Vert _{H}^{2}}.  \label{62-0}
\end{equation}%
We note that by the Rayleigh-Ritz theorem (see e.g., \cite{Reed-Simon}) the
first eigenvalue of the operator $B_{2}B_{2}^{\ast }$ is given by 
\begin{equation}
\lambda _{1}^{B}=\inf_{y\in H,\mbox{ }y\neq 0}\frac{(B_{2}B_{2}^{\ast
}y,y)_{H}}{\left\Vert y\right\Vert _{H}^{2}}=\inf_{y\in H,\mbox{ }y\neq 0}%
\frac{\left\Vert B_{2}^{\ast }y\right\Vert _{H}^{2}}{\left\Vert y\right\Vert
_{H}^{2}}  \label{63}
\end{equation}%
which is positive due to (\ref{64-0}) since this means that $%
(B_{2}B_{2}^{\ast }y,y)_{H}=\left\Vert B_{2}^{\ast }y\right\Vert _{H}^{2}>0$
for all $y\in H,$ $y\neq 0$. Moreover, (\ref{64}) implies (\ref{62-0}) and
so, we get $\Gamma \geq 0$ as claimed.\hfill $\square $

As examples of such kernels $k$ we can mention the Gaussian kernel $k(x,\xi
)=\exp \left( -\frac{||x-\xi ||^{2}}{2\sigma ^{2}}\right) ,$ $\sigma >0,$ or
the Newtonian kernel $k(x,\xi )=\alpha \ln |x-\xi |$ (e.g., for $N=2),$ $%
\alpha >0.$

In conclusion, Theorem \ref{T2.1} can be applied to obtain the existence and
uniqueness of the solution to the Riccati equation (\ref{6}) corresponding
to system (\ref{21-0})-(\ref{21-3}).

\end{document}